\documentclass[amscd, amssymb,12pt]{article}
\usepackage{amsmath, amsopn,amsthm,amsfonts}
\normalbaselines
\newtheorem{Theorem}{Theorem}[section]
\newtheorem{Definition}[Theorem]{Definition}
\newtheorem{Proposition}[Theorem]{Proposition}
\newtheorem{Lemma}[Theorem]{Lemma}
\newtheorem{Corollary}[Theorem]{Corollary}

\newtheorem{Remark}[Theorem]{Remark}
\newtheorem{Example}[Theorem]{Example}

\def\A{{\cal A}}
\def\B{{\cal B}}
\def\cD{{\cal D}}
\def\H{{\cal H}}
\def\L{{\cal L}}
\def\O{{\cal O}}
\def\P{{\mathbb P}}

\def\ca{c\`adl\`ag}
\def\K{{\cal K}}
\newcommand{\esp}{{\rm I\!E}}
\newcommand{\myref}[1]{(\ref {#1})}

\newcommand{\I}{{\bf 1}}
\newcommand{\R}{{\rm I\!R}}
\newcommand{\espe}{{\rm I\!E}}
\newcommand{\der}[2]{{\frac{ \textstyle \partial #1}{ \textstyle \partial #2}}}

\newcommand{\Tr}{\mathop{\rm Tr}}
\def\qed{\hbox{\hskip 6pt\vrule width6pt height7pt depth1pt \hskip1pt}}
\newcommand{\ioi}{\int_0^{+\infty}}
\def\<{\left\langle} \def\>{\right\rangle}
\let\D\Delta
\let\d\delta
\let\g\gamma
\let\k\kappa
\let\la\lambda

\let\wh\widehat
\let\er\eqref
\def\R{\mathbb R}
\def\Qp{Q^{1/2}}
\def\fp{\frac12}
\def\Im{\mathop{\rm Im}}

\def\beq{\begin{equation}}
\def\bmlg{\begin{multline*}}
\def\mx#1{\matrix#1\endmatrix}
\def\pmx#1{\left(\mx{#1}\right)}
\def\ca{c\`adl\`ag}
\def\<{\left\langle} \def\>{\right\rangle}
\let\D\Delta
\let\d\delta
\let\g\gamma
\let\k\kappa
\let\la\lambda
\def\R{\mathbb R}
\def\Qp{Q^{1/2}}
\def\fp{\frac12}
\def\Im{\mathop{\rm Im}}

\def\beq{\begin{equation}}
\def\mx#1{\begin{matrix}#1\end{matrix}}
\def\pmx#1{\left(\mx{#1}\right)}

\begin{document}
\baselineskip=17pt

\title{{\bf Large deviations for stochastic PDE with L\'evy noise
\footnote{ Supported by the  MNiSW project 1PO 3A 034 29
``Stochastic evolution equations with L\'evy noise'', EC FP6 Marie
Curie ToK programme SPADE2, and NSF grants DMS-0500270 and DMS-0856485.}}}
\author{
    \textsc{Andrzej \'Swi\c{e}ch}\footnote{ Corresponding author. E-mail:
swiech@math.gatech.edu.}\\
    \textit{School of Mathematics, Georgia Institute of Technology}\\
\textit{
Atlanta, GA 30332, U.S.A.} \vspace{.2cm}\\
    {\small AND} \vspace{.2cm} \\
    \textsc{Jerzy Zabczyk}\\
    \textit{Institute of Mathematics, Polish Academy of Sciences}\\
\textit{Sniadeckich 8, 00-950 Warsaw, Poland}
    }
\date{}
\maketitle

\begin{abstract}
We prove a large deviation principle result for solutions of abstract stochastic evolution equations perturbed by small L\'evy noise. We use general
large deviations theorems of Varadhan and Bryc, viscosity solutions of integro-partial differential equations in Hilbert spaces, and deterministic
optimal control methods. The Laplace limit is identified as a viscosity solution of a Hamilton-Jacobi-Bellman equation of an associated
control problem. We also establish exponential moment estimates for solutions of stochastic evolution equations driven by L\'evy noise. General results
are applied to stochastic hyperbolic equations perturbed by subordinated Wiener process.
\end{abstract}

\vspace{.2cm}

\vspace{.2cm}

\section{Introduction}

\label{INTRO}
\setcounter{equation}{0}

Let $L(t)$ be a square integrable L\'evy martingale on a Hilbert space $H$, starting from $0$, defined on a complete probability space $(\Omega,{\cal
F},{\mathbb P})$ with a normal filtration ${\cal F}_t$. It is well known, see e.g. \cite{PZ}, that
\begin{equation}\label{L1}
L(t)=\int_0^t\int_H z\hat\pi(ds,dz) + W(t)
\end{equation}
where $W$ is an $H$ valued Wiener process,  independent of the
compensated random measure $\hat\pi(ds,dz)=\pi(ds,dz)-ds\nu(dz)$
with the intensity measure $\nu$, satisfying
$$
\int_{H} \|z\|^2 \nu (dz) < +\infty .
$$
Here
 $$
 \pi(]0,t], \Gamma) = \#\{s\in ]0,t]; L(t)-L(t-) \in \Gamma\},
 $$
is the random measure of jumps of the process $L$, see e.g. \cite{S}, \cite{B} and \cite{PZ}. Define
$$
L_n(t)=\frac{1}{n}L(nt),
$$
and note that
$$
\esp \|L_{n}(t)\|^{2}= {\frac{t}{n}} \int_{H}\|z\|^{2}\nu(dz).
$$

We study large deviation principle for the family of processes
$\{X_n\}$ satisfying
\begin{equation}
\label{a1} dX_n(s)=(-AX_n(s)+F(X_n(s)))ds+G(X_n(s-))dL_n(s),\quad
X_n(0)=x\in H,
\end{equation}
where $A$ is a linear, densely defined, maximal monotone operator
in $H$ and $F, G$ are certain continuous functions. These abstract
stochastic differential equations may be for instance semilinear
stochastic PDE with small L\'evy noise. For the theory of such
equations we refer to \cite{PZ} and the references therein. We
excluded from our considerations the Gaussian part of the noise.
If $L$ is a Wiener process, large deviation results are well
known, see e.g. \cite{cero, chmil, pchow, DPZ, FK, frei, kalxi,
Pe, so, SriSun, Sw} and the references therein. We think that our
methods, combined with the techniques of \cite{Sw}, should apply
to the general case, however we do not attempt to do it here.
Thus, we will always assume that
\begin{equation}
\label{levy} L_n(t)=\frac{1}{n}L(nt),\,\,{\rm where}\,\, L(t)=\int_0^t\int_H z\hat\pi(ds,dz).
\end{equation}

There are two types of large deviation results; at a single time,
i.e. for $X_n(T)$ with $T$ fixed, and in the path space, i.e. for
$X_n(\cdot)$. Our goal is to show the large deviation principle
and identify the rate function for the single time case since this
is where the PDE theory is used. Once this is done a general
strategy to pass to the path space case can be found in \cite{FK}.
Such a strategy was employed in \cite{Sw} when $L$ was a Wiener
process. We don't know if it can be successful here.

The problem of large deviations for infinite dimensional processes
with jumps seems to be wide open although for the finite
dimensional spaces basic results are presented in \cite{W}. We are
only aware of three papers that specifically address it in the path
space. In \cite{A1} the large deviation principle is proved for a family of
Banach space valued L\'evy processes and in \cite{RZ} for
solutions of linear evolution equations of type (\ref{a1}) with
additive L\'evy noise and the operator $A$ with a discrete
spectrum. Paper \cite{XZ} deals with the case of two-dimensional
stochastic Navier-Stokes equations driven by additive L\'evy noise.
We also refer to \cite{A2, FK} for related results.

Our approach uses the classical theorems of Varadhan and Bryc
\cite{DPE}. According to them the processes $X_n$ satisfy the
large deviation principle in a metric space $E$ if and only if the
family $\{X_n\}$ is exponentially tight and the {\it Laplace
limit}
\[
\Lambda(g)=\lim_{n\to\infty}\frac{1}{n}\log \esp e^{ng(X_n)}
\]
exists for all $g\in C_b(E)$. We will choose $E$ to be any Hilbert space $V$ such that $H\subset V$ and $H\hookrightarrow V$ is compact. Our
main result, the existence of the Laplace limit and its identification, will be a consequence of a much more general result about convergence of
viscosity solutions of certain integro-PDE in $H$ to the viscosity solution of the limiting first order Hamilton-Jacobi-Bellman (HJB)
equation.\vspace{2mm}

After recalling basic definitions and introducing main hypotheses in Section \ref{prelim},   exponential estimates and continuous dependence estimates
for solutions of (\ref{a1}) are established in Section \ref{estspdelevy}, see
Proposition \ref{msol} and Proposition \ref{apt1}. In the proofs we use  a new
result on
convergence of solutions of  the equation (\ref{a1}) with $A$ replaced by Yosida approximations of $A$. Associated nonlinear PDE in Hilbert
spaces are investigated in Section \ref{intpde}. The fact that functions
$$
 v_n(t,x) = \frac{1}{n}\log \esp e^{ng(X_n(T))},
$$
where $X_n$ solves (\ref{a1a}), are viscosity solutions of proper
nonlinear PDE, is the content of Theorem \ref{integropde}.
Moreover Theorem \ref{vercontpr} establishes existence of a
viscosity solution to the limiting HJB equation. The main results
on the Laplace limits  are subjects  of Theorem \ref{thcompcp1},
Theorem \ref{thlimit}, and Corollary \ref{corlimit} of Section
\ref{bcontcomparison}. Finally Theorem \ref{thldp} states
conditions under which the large deviation principle holds for
solutions of (\ref{a1}). Various examples are discussed in
Sections \ref {examplesnoise} and \ref{wave}. In the Appendix we
give a proof of the  convergence result used in Section
\ref{estspdelevy}.

\section{Preliminaries}

\label{prelim}
\setcounter{equation}{0}

\subsection{Basic definitions and assumptions}
Throughout this paper $H$ will be a real separable Hilbert space equipped with the inner product
$\langle \cdot,\cdot\rangle$ and the norm $\|\cdot\|$.
We recall that $A$ is a linear, densely defined, maximal monotone operator in $H$.

Let $B$ be a bounded, linear, positive, self-adjoint operator on
$H$ such that $A^*B$ is bounded on $H$ and
\begin{equation}
\label{bcond} \langle (A^* B + c_0 B)x,x \rangle \geq 0 \;\;\;\;\;\;
\text{for all} \,\, x\in H
\end{equation}
for some $c_0\geq 0$. Such an operator always exists, for instance
$B=((A+I)(A^*+I))^{-1/2}$ (see \cite{Re}). We refer to \cite{CL4}
for various examples of $B$. Using the operator $B$ we define for
$\gamma>0$ the space $H_{-\gamma}$ to be the completion of $H$
under the norm
\[
\|x\|_{-\gamma}=\|B^{\frac{\gamma}{2}}x\|.
\]
Let $\Omega\subset [0,T]\times H$. We say that $u:\Omega\to
\mathbb{R}$ is $B$-upper-semicontinuous (respectively,
$B$-lower-semicontinuous) on $\Omega$ if whenever $t_n\to t$,
$x_n\rightharpoonup x$, $Bx_n\to Bx$, $(t,x)\in \Omega$, then
$\limsup_{n\to+\infty}u(t_n,x_n)\leq u(t,x)$ (respectively,
$\liminf_{n\to+\infty}u(t_n,x_n)\geq u(t,x)$). The function $u$ is
$B$-continuous on $\Omega$ if it is $B$-upper-semicontinuous and
$B$-lower-semicontinuous on $\Omega$.

The following assumptions will be made about the functions $F:H\to
H$ and $G:H\to L(H)$, where $L(H)$ is the space of bounded linear
operators on $H$:

\begin{equation}\label{as1}
\|F(0)\|\leq M,\quad\|F(x)-F(y)\|\leq M\|x-y\|_{-1} \quad\mbox{for
all}\,\,x,y\in H,
\end{equation}
\begin{equation}\label{as2}
\|G(x)-G(y)\|\leq M\|x-y\|_{-1}\quad\mbox{for all}\,\,x,y\in H,
\end{equation}
\begin{equation}\label{as3}
\|G(x)\|\leq M\quad\mbox{for all}\,\,x\in H
\end{equation}
for some $M\geq 0$, and
\begin{equation}\label{as4}
\int_H \|z\|^2e^{K\|z\|}\nu(dz)<+\infty\quad\mbox{for
every}\,\,K>0.
\end{equation}
Condition (\ref{as4}) is equivalent to the requirement that the
noise process has exponential moments:
$$
\esp e^{K\|L(t)\|} < +\infty,\quad\mbox{for all}\,\,t,K>0.
$$

If \eqref{as4} holds then the Laplace transform of the process $L$ is well defined. Namely if $L$ is given by \eqref{L1} and $Q_{W}$ is the covariance of
$W$, then
\begin{equation}
\esp e^{\langle p, L(t) \rangle} = e^{t H(p)} \,\,{\rm where}\,\, H(p)= 1/2 \langle Q_{W}p, p \rangle +\int_H\left[e^{\langle p,z\rangle}-1-\langle
p,z\rangle\right]\nu(dz),\,\,p\in H.
\end{equation}
We set
$$
H_{0}(p) =\int_H\left[e^{\langle p,z\rangle}-1-\langle p,z\rangle\right]\nu(dz),\,\,p\in H,
$$
if $L$ is without the Gaussian part as in \eqref{levy}.
\begin{Remark}
\label{remBstrong} If instead of (\ref{bcond}) we suppose that
\begin{equation}
\label{bcondstrong} \langle (A^* B + c_0 B)x,x \rangle \geq
\|x\|^2 \;\;\;\;\;\; \text{for all} \,\, x\in H
\end{equation}
then (\ref{as1}) can be replaced by a weaker condition
\begin{equation}\label{as1strong}
\|F(0)\|\leq M,\quad\|F(x)-F(y)\|\leq M\|x-y\| \quad\mbox{for
all}\,\,x,y\in H.
\end{equation}
We refer the reader to \cite{CL4} for examples of operators
satisfying (\ref{bcondstrong}) and to \cite{Re} for conditions
guaranteeing the existence of $B$ for which (\ref{bcondstrong})
holds.
\end{Remark}

We will need the following simple fact which we record for future
use.

\begin{Lemma}\label{intform}
If $f\in C^2(H)$ then for every $x,y\in H$
\[
f(x+y)=f(x)+\langle Df(x),y\rangle+\int_0^1\int_0^1\langle D^2f(x+s\sigma y)y,y\rangle\sigma ds d\sigma.
\]
\end{Lemma}

For a square integrable martingale $M$ we will denote by $\langle
M,M\rangle_t$ its angle bracket and by $[M,M]_t$ its quadratic
variation (see \cite{P}, p.~57, or \cite{M}, p. 150). It is easy to
see that $\langle L(nt),L(nt)\rangle_t=cnt$ for some $c>0$.

For a Hilbert space $Z$ we will be using the following function spaces.
\[
C_b(Z)=\{u:Z\to{\mathbb R}: u\,\,\mbox{is continuous and
bounded}\},
\]
\[
{\rm Lip}_b(Z)=\{u\in C_b(Z): u\,\,\mbox{is Lipschitz continuous}\},
\]
\[
C^{2}(Z)=\{u:Z\to{\mathbb R}: Du,D^2u\,\,\hbox{are continuous}\},
\]
\[
C^{1,2}((0,T)\times Z)=\{u:(0,T)\times Z\to{\mathbb R}:
u_t,Du,D^2u\,\,\hbox{are continuous}\},
\]
\[
C^{2}_{\rm uc}(Z)=\{u:Z\to{\mathbb R}: u,Du,D^2u\,\,\hbox{are
uniformly continuous}\},
\]
where $Du, D^2u$ denote the Fr\'echet derivatives of $u$ with
respect to the spatial variable.

We will denote by $S(\cdot)$ the $C_0$-semigroup generated by $-A$. For $\lambda>0$ we denote by $A_\lambda$ the {\it Yosida approximation} of $A$,
$A_\lambda=\lambda A R_\lambda$, where $R_\lambda=(\lambda I+A)^{-1}$. The $C_0$-semigroup generated by $-A_\lambda$ will be denoted by
$S_\lambda(\cdot)$. Both $S(\cdot)$ and $S_\lambda(\cdot)$ are semigroups of contractions. It is well known (see for instance \cite{Pazy}) that
\begin{equation}
\label{resolv} \|R_\lambda\|\leq \frac{1}{\lambda},\quad
\hbox{and}\,\,\,\lim_{\lambda\to+\infty}\lambda R_\lambda
x=x\,\,\,\hbox{for}\,\,x\in H.
\end{equation}
For $C\in L(H)$ we will denote by $\| C \|_{\rm HS}$ its
Hilbert-Schmidt norm.

\subsection{Viscosity solutions}

\label{VS}
\setcounter{equation}{0}

To minimize the technicalities we will be using a slightly
simplified definition of viscosity solution. This simplified
definition will be enough since in this paper we only deal with
bounded solutions. We also point out that Definition
\ref{defvisc2-t} applies to terminal value problems.

\begin{Definition}
\label{dftestftd} A function $\psi $ is a test function if $\psi
=\varphi + h(\|x\|)$, where:

\begin{enumerate}
\item[(i)] $\varphi \in C^{1,2}\left( \left( 0,T\right) \times
H\right)$, is $B$-lower semicontinuous, $\varphi,\varphi_t,
D\varphi, D^2\varphi$, $A^\ast D\varphi$ are uniformly continuous
on $[\epsilon,T-\epsilon]\times H$ for every $\epsilon>0$, and
$\varphi$ is bounded on every set $[\epsilon,T-\epsilon]\times
\{\|x\|_{-1}\le r\}$.

\item[(ii)] $h\in C^2([0,+\infty))$ is such that $h'(0)=0,
h'(r)\ge 0$ for $r\in (0,+\infty)$, and $h,h',h''$ are uniformly continuous on
$[0,+\infty)$.
\end{enumerate}

\end{Definition}

We will be concerned with terminal value problems for integro-PDE of the form
\begin{equation}
v_t-\langle Ax,Dv\rangle+\,F(t,x,Dv,v(t,\cdot)) =0\quad\hbox{in}\,\,
(0,T)\times H, \label{eq:CP}
\end{equation}
where $F:(0,T)\times H\times H\times C^2_{\rm uc}(H) \to{\mathbb
R}$.

\begin{Definition}
\label{defvisc2-t} A locally bounded $B$-upper semicontinuous
function $u:(0,T)\times H\to {\mathbb R}$ is a viscosity subsolution of
(\ref{eq:CP}) if whenever $u -\varphi-h(\|\cdot\|)$ has a maximum over $(0,T)\times H$
at a
point $(t,x)$ for some test functions $\varphi, h(\|y\|)$ then
\[
\psi_t(t,x) -\langle x,A^*D\varphi (t,x)\rangle
+F(t,x,D\psi(t,x),\psi(t,\cdot))\geq 0,
\]
where $\psi(s,y) =\varphi(s,y)+h(\|y\|).$

A locally bounded $B$-lower semicontinuous function $u:(0,T)\times
H\to {\mathbb R}$ is a viscosity supersolution of (\ref{eq:CP}) if whenever
$u +\varphi+h(\|\cdot\|)$ has a minimum over $(0,T)\times H$ at a point $(t,x)$ for
some test
functions $\varphi, h(\|y\|)$ then
\[
\psi_t(t,x)
+\langle x,A^*D\varphi (t,x)\rangle
+F(t,x,D\psi(t,x),\psi(t,\cdot))\leq 0,
\]
where $\psi(s,y) =-\varphi(s,y)-h(\|y\|).$

A viscosity solution of (\ref{eq:CP}) is a function which is both
a viscosity subsolution
and a viscosity supersolution.
\end{Definition}

\section{Estimates for solutions of stochastic PDE with\\ L\'evy noise}

\label{estspdelevy}

\setcounter{equation}{0}

In this section we recall basic facts and show various estimates
about mild solutions of the equations,
\begin{equation}
\label{a1a} dX_n(s)=(-AX_n(s)+F(X_n(s)))ds+G(X_n(s-))dL_n(s),\quad
X_n(t)=x\in H,
\end{equation}
on a fixed time interval $[0,T]$, where $L_n$ are the processes defined in (\ref{levy}).\vspace{2mm}

\noindent Let us recall that if \eqref{levy} holds then
\begin{equation}
\esp e^{\langle p, L_{n}(t) \rangle} = e^{nt H_{0}({\frac{p}{n}})} = e^{nt\int_H\left[e^{{\frac{1}{n}}\langle p,z\rangle}-1-{\frac{1}{n}}\langle
p,z\rangle\right]\nu(dz)},\,\,p\in H.
\end{equation}
The covariance operator of the process $L$ will be denoted by $Q$ and then the covariance operator of $L_n$ is $\frac{1}{n} Q$.\vspace{2mm}

\noindent We refer the readers to Chapter 9 of \cite{PZ} for the definition of a mild solution. We will also need solutions $X_n^m$ of the equations
\begin{equation}
\label{m1a} dX_n^m(s)=(-A_mX_n^m(s)+F(X_n^m(s)))ds+G(X_n^m(s-))dL_n(s),\quad X_n^m(t)=x\in H,
\end{equation}
where the operators $A_m$  are Yosida approximations of $A$ for $\lambda=m= 1,2,\ldots$.
\begin{Proposition}\label{msol}
Let $0\le t\le T$. Let (\ref{as4}) be satisfied and let
\begin{equation}\label{cond24}
\|G(x)-G(y)\|,\|F(x)-F(y)\|\leq C\|x-y\| \quad\mbox{for
all}\,\,x,y\in H,
\end{equation}
for some $C\geq 0$. Then:

(i) There exists a unique mild solution $X_n$ of (\ref{a1a}). The
solution $X_n$ has a \ca\ modification.

(ii) If $X_n^m$ is the solution of (\ref{m1a}) then
\begin{equation}\label{yapr}
\lim_{m\to+\infty}\esp\left(\sup_{t\leq s\leq
T}\|X_n^m(s)-X_n(s)\|^2\right)=0.
\end{equation}

(iii) If in addition (\ref{as3}) holds then there exist constants
$c_1>0, c_2>0$ (depending only on $T,M$, with $c_2$ depending also
on $\|x\|$) such that
\begin{equation}\label{eme1}
\esp\left(\sup_{t\leq s\leq T} e^{nc_1\|X_n(s)\|}\right)\leq
e^{nc_2}.
\end{equation}
\end{Proposition}

\begin{Remark}\label{remmom}
It follows from the proof that (\ref{eme1}) is also satisfied for
the processes $X_n^m$ with the same constants $c_1,c_2$. In
particular this implies that there exists a constant $C(\|x\|,T)$
such that for every $n,m$
\begin{equation}\label{momest}
\esp\left(\sup_{t\leq s\leq T} e^{c_1\|X_n^m(s)\|}\right)\leq
C(\|x\|,T)
\end{equation}
with the same estimate being also true for the processes $X_n$.
\end{Remark}
{\bf Proof}. $(i)$ This is a standard result, see Theorem 9.29 in
\cite{PZ}.

$(ii)$ We will need two general results on convergence of stochastic and deterministic convolutions, Propositions \ref{apt1} and
\ref{apt2}. The proof of Proposition \ref{apt1} will be postponed to the Appendix and the classical proof of Proposition \ref{apt2} will be omitted.

Denote by  $\cal L$ the space of all predictable processes $\psi (\cdot)$ whose values are linear operators from the space $Q^{1/2}(H)$ into $H$,
equipped with the scalar product
$$
<\psi_1, \psi_2>_{\cal L} = \sum_{n=1}^{+\infty}<\psi_{1}Q^{1/2}e_n ,\psi_{2}Q^{1/2}e_n  >_{H},\,\,\,\psi_1, \psi_2 \in {\cal L}.
$$
Here $(e_n)$ is any orthonormal basis in $H$. Moreover two operators on $H$, even unbounded,  identical on $Q^{1/2}(H)$, are identified. The norm on
$\cal L$ is given by the formula.
\[
|\psi|_{1}= \left(\esp \int_0^T \bigl\|\psi(s)\Qp\bigr\|^2_{\rm HS}\,ds\right)^{1/2} <+\infty .
\]
\begin{Proposition}\label{apt1}
Let $L(t)$ be a square integrable L\'evy martingale in $H$ with the covariance operator $Q$, and $\psi \in \cal L$. Then the processes
\begin{equation}\label{ap11125}
\int_0^t S(t-s)\psi(s)\,dL(s),\,\, \int_0^t S_\la(t-s)\psi(s)\,dL(s),\,\, t\in[0,T],\,\,\,\,\,\la> 0,
\end{equation}
have \ca\ modifications and
\begin{equation}\label{ap1}
\lim_{\la\to+\infty} \esp \sup_{0\le t\le T} \big\|\int_0^t
S(t-s)\psi(s)\,dL(s) - \int_0^t
S_\la(t-s)\psi(s)\,dL(s)\big\|^2=0.
\end{equation}
\end{Proposition}

\begin{Proposition}\label{apt2}
Assume that $\psi$ is an $H$-valued predictable process such that
$$
\esp \int_0^T \|\psi(s) \|^2 \,ds<+\infty
$$
Then the processes
\[
\int_0^t S(t-s)\psi(s)\,ds,\,\,\, \int_0^t
S_\la(t-s)\psi(s)\,ds,\quad
  t\in[0,T],\,\,\,\la> 0,
\]
have continuous modifications and
\[
\lim_{\la\to+\infty} \esp \sup_{0\le t\le T} \big\|\int_0^t
S(t-s)\psi(s)\,ds - \int_0^t S_\la(t-s)\psi(s)\,ds\big\|^2=0.
\]
\end{Proposition}

\medskip
We can now proceed with the proof of $(ii)$. Let $\cal X$ denote the space of all \ca, adapted to the filtration ${\cal F}_t$, $H$-valued processes $X$,
equipped with the norm $|\cdot|_0$:
\[
|X|_0 = \left(\esp\sup_{t\leq T}\|X(t)\|^{2}\right)^{1/2} .
\]
Define transformations ${\cal K}_{n},\,\,$ ${\cal K}_{nm},\,\,$ ,\,\,$n,m =1,2,\ldots$ by the formulae,
\[
{\cal K}_{n}(X)(t) = S(t)X_0+\int_0^t S(t-s)F(X(s))ds +\int_0^tS(t-s)G(X(s-))dL_n(s),
\]
\[
{\cal K}_{nm} (X)(t) = S_m(t)X_0+\int_0^t S_m(t-s)F(X(s))ds +\int_0^tS_m(t-s)G(X(s-))dL_n(s).
\]
It will follow from the first part of the proof of Proposition \ref{apt1} that the processes ${\cal K}_{n}(X)$, ${\cal K}_{nm}(X)$ have \ca\ modifications.
Moreover, as in the proof of existence of mild solutions, see e.g. \cite {PZ} and using arguments similar to the proof of (\ref{ap3})  one can show that
for arbitrary $\alpha \in (0,1)$ there exists $T_{\alpha}$ such that all transformations ${\cal K}_{n},\,\,{\cal K}_{n}$ satisfy Lipschitz conditions on
$\cal X$ with a constant smaller than $\alpha$. Moreover processes $X_n ,$\,$X_n^m ,$ are unique solutions in $\cal X$ of the following fixed point
problems
$$
X={\cal K}_{n}(X),\,\,\,\,\,X={\cal K}_{mn}(X).
$$
Therefore, it is easy to see, that to prove the results it is enough to show that for each $X\in \cal X$,
$$
\lim_{m}{\cal K}_{mn}(X)= {\cal K}_{n}(X),
$$
and this follows from Proposition \ref{apt1},  \ref{apt2}. The case of arbitrary $T>0$
follows by repeating the same argument on intervals $[0,
T_{\alpha}]$,\,\,$[T_{\alpha}, 2T_{\alpha}]$,\ldots,$[(k-1)T_{\alpha}, kT_{\alpha}]$,
where $kT_{\alpha} >T$.

$(iii)$ Without loss of generality we will assume that $t=0$. We
will denote by $\pi_n(dt,dz)$, respectively $\pi_n^k(dt,dz), k\geq
1$, the Poisson random measure for the process $L(nt)$,
respectively $L^k(nt)$, where $L^k(nt)$ is the process $L(nt)$
with jumps restricted to size $k$. It is easy to see that the
intensity measure of $L(nt)$ is equal to $n\nu(dz)$ and the
intensity measure of $L^k(nt)$ is equal to $n\nu^k(dz)$, where
$\nu^k(dz)=\chi_{\{\|z\|\leq k\}}\nu(dz)$.

Denote by $X_n^{mk}, m,k=1,2,...$ the solution of (\ref{a1a}) with
$A$ replaced by $A_m$ and $L_n$ replaced by $L_n^k$, where
$L_n^k=\frac{1}{n}L^k(nt)$. We will show (\ref{eme1}) for the
processes $X_n^{mk}$ and then pass to the limit as $k\to+\infty$
and $m\to+\infty$.

Let $h:{\mathbb R}\to {\mathbb R}$ be a smooth even function such
that $h(0)=1$, $h$ is increasing on $(0,+\infty)$, $h'(0)=0$,
$|h'(r)|\leq 1, h(r)\geq (1+r)/2$ for $r>0$. (We can take for
instance $h(r)=\sqrt{1+r^2}$.) For $l>0$ denote by $\tau_l$ the
exit time of $X_n^{mk}$ from $\{\|y\|\leq l\}$. Let $\alpha>0$ be
a number which will be specified later. By Ito's formula, see
\cite{M}, Theorem 27.2, p. 190, we have
\begin{equation}\label{eme2}
\begin{split}
&e^{ne^{-\alpha (s\wedge\tau_l)}h(\|
X_n^{mk}(s\wedge\tau_l)\|)}=e^{n h(\|x\|)}
\\
& \quad -\int_0^{s\wedge\tau_l}\alpha n e^{-\alpha
r}h(\|X_n^{mk}(r)\|) e^{ne^{-\alpha r}h(\|X_n^{mk}(r)\|)}dr
\\
& \quad+\int_0^{s\wedge\tau_l}ne^{-\alpha r} e^{ne^{-\alpha
r}h(\|X_n^{mk}(r)\|)}h'(X_n^{mk}(r))\langle-A_mX_n^{mk}(r)+F(X_n^{mk}(r)),
\frac{X_n^{mk}(r)}{\|X_n^{mk}(r)\|}\rangle dr
\\
& \quad + \int_0^{s}n e^{-\alpha r} e^{n e^{-\alpha
r}h(\|X_n^{mk}(r-)\|)}h'(X_n^{mk}(r-)){\bf 1}_{[0,\tau_l]}\langle
\frac{X_n^{mk}(r-)}{\|X_n^{mk}(r-)\|},
G(X_n^{mk}(r-))dL_n^k(r)\rangle
\\
& \quad + \int_0^{s}\int_H{\bf 1}_{[0,\tau_l]} \bigg[e^{n
e^{-\alpha r}h(\|X_n^{mk}(r-)+\frac{1}{n}G(X_n^{mk}(r-))z\|)}-
e^{n e^{-\alpha r}h(\|X_n^{mk}(r-)\|)}
\\
& \quad\quad -e^{-\alpha r}e^{n e^{-\alpha
r}h(\|X_n^{mk}(r-)\|)}h'(X_n^{mk}(r-))\langle
\frac{X_n^{mk}(r-)}{\|X_n^{mk}(r-)\|},
G(X_n^{mk}(r-))z\rangle\bigg]\pi_n^k(dr,dz).
\end{split}\end{equation}
To proceed further we  compensate the measure $\pi$ and recall that stochastic integrals
with respect to the compensated random measures form martingales. Thus taking expectation
in (\ref{eme2}), using (\ref{as3}),
(\ref{cond24}), martingale property, the fact that $\langle -A_m
y,y\rangle\leq 0$ for $y\in H$ and $1+r\leq 2h(r)$, we therefore
obtain
\begin{equation}\label{eme221}
\begin{split}
&\esp e^{ne^{-\alpha (s\wedge\tau_l)}h(\|
X_n^{mk}(s\wedge\tau_l)\|)}\leq e^{n h(\|x\|)}
\\
& \quad\quad +\esp\int_0^{s\wedge\tau_l}n e^{-\alpha r}
e^{ne^{-\alpha
r}h(\|X_n^{mk}(r)\|)}\left[C(1+\|X_n^{mk}(r)\|)-\alpha
h(\|X_n^{mk}(r)\|)\right] dr
\\
& \quad\quad + \esp\int_0^{s\wedge\tau_l}\int_H n\bigg|e^{n
e^{-\alpha r}h(\|X_n^{mk}(r)+\frac{1}{n}G(X_n^{mk}(r))z\|)}- e^{n
e^{-\alpha r}h(\|X_n^{mk}(r)\|)}
\\
& \quad\quad\quad -e^{-\alpha r}e^{n e^{-\alpha
r}h(\|X_n^{mk}(r)\|)}h'(X_n^{mk}(r))\langle
\frac{X_n^{mk}(r)}{\|X_n^{mk}(r)\|},
G(X_n^{mk}(r))z\rangle\bigg|\nu(dz)dr
\\
& \quad\leq e^{n h(\|x\|)} +\esp\int_0^{s\wedge\tau_l}n e^{-\alpha
r} e^{ne^{-\alpha r}h(\|X_n^{mk}(r)\|)}(2C-\alpha)
h(\|X_n^{mk}(r)\|) dr
\\
& \quad\quad\quad\quad\quad\,\,\,+\esp\int_0^{s\wedge\tau_l}I(r)
dr,
\end{split}\end{equation}
where $I(r)$ is the integrand of the last term in the middle line
of (\ref{eme221}). Applying Lemma \ref{intform} to the function
$f(x)=e^{ne^{-\alpha r}h(\|x\|)}$ we have
\begin{equation}\label{formi}
\begin{split}
& I(r)=\int_H\bigg|\int_0^1\int_0^1 n\langle
D^2f\left(X_n^{mk}(r)+\frac{t\sigma}{n}G(X_n^{mk}(r))z\right)\frac{1}{n}G(X_n^{mk}(r))z,
\\
& \quad\quad\quad\quad\quad\quad
\frac{1}{n}G(X_n^{mk}(r))z\rangle\,\sigma dt\,
d\sigma\bigg|\nu(dz).
\end{split}
\end{equation}
Elementary calculation gives us
\[
D^2f(x)=ne^{-\alpha r}e^{ne^{-\alpha
r}h(\|x\|)}(n\psi_1(x)+\psi_2(x)),
\]
where
\[
\psi_1(x)=e^{-\alpha
r}(h'(\|x\|))^2\frac{x}{\|x\|}\otimes\frac{x}{\|x\|},
\]
\[
\psi_2(x)=\left(h''(\|x\|)-\frac{h'(\|x\|)}{\|x\|}\right)\frac{x}{\|x\|}\otimes\frac{x}{\|x\|}
+\frac{h'(\|x\|)}{\|x\|}I.
\]
We observe that both $\psi_1,\psi_2$ are bounded as functions from
$H$ to $L(H)$. Therefore
\begin{equation}\label{formi1}
\begin{split}
& I(r)\leq e^{ne^{-\alpha r}h(\|X_n^{mk}(r)\|)}
\int_H\int_0^1\int_0^1 M^2e^{-\alpha r}
\\
& \quad\quad\quad\quad e^{n e^{-\alpha
r}|h(\|X_n^{mk}(r)+\frac{t\sigma}{n}G(X_n^{mk}(r))z\|)-h(\|X_n^{mk}(r)\|)|}(n\|\psi_1\|_\infty
+\|\psi_2\|_\infty)\|z\|^2dt\,d\sigma\nu(dz)
\\
& \quad\quad\leq e^{ne^{-\alpha r}h(\|X_n^{mk}(r)\|)} \int_H
M^2e^{-\alpha r}e^{M\|z\|}(n\|\psi_1\|_\infty
+\|\psi_2\|_\infty)\|z\|^2\nu(dz)
\\
& \quad\quad\leq nM_1e^{-\alpha r}e^{ne^{-\alpha
r}h(\|X_n^{mk}(r)\|)}\int_H\|z\|^2e^{M\|z\|}\nu(dz) \leq
nM_2e^{-\alpha r}e^{ne^{-\alpha r}h(\|X_n^{mk}(r)\|)}
\end{split}
\end{equation}
for some $M_1,M_2>0$. Plugging (\ref{formi1}) into (\ref{eme221}),
choosing $\alpha=2C+M_2+1$ and recalling that $h(r)\geq 1$ we thus
obtain
\[
\esp e^{ne^{-\alpha (s\wedge\tau_l)}h(\|
X_n^{mk}(s\wedge\tau_l)\|)}+\esp\int_0^{s\wedge\tau_l}ne^{-\alpha
r}e^{ne^{-\alpha r}h(\|X_n^{mk}(r)\|)}dr\leq e^{n h(\|x\|)}
\]
which in particular implies that
\[
\esp e^{ne^{-\alpha s}h(\| X_n^{mk}(s\wedge\tau_l)\|)}\leq e^{n
h(\|x\|)}.
\]
Since $\lim_{l\to+\infty}(T\wedge\tau_l)=T$ a.s., letting
$l\to+\infty$ and using Fatou's lemma  we obtain
\[
\esp e^{ne^{-\alpha s}h(\| X_n^{mk}(s)\|)}\leq e^{n h(\|x\|)}.
\]
We can now send $k\to+\infty$, employ once again Fatou's lemma and
the fact that $X_n^{mk}(s)\to X_n^{m}(s)$ a.s. (at least along a
subsequence). This can be shown using the arguments from the proof
of $(ii)$. This way we arrive at
\begin{equation}\label{expest1}
\esp e^{ne^{-\alpha s}h(\| X_n^{m}(s)\|)}\leq e^{n h(\|x\|)}.
\end{equation}

We can now go back to Ito's formula (\ref{eme2}) but apply it to
the function $e^{\frac{n}{2}e^{-\alpha r}h(\|x\|)}$, the process
$X_n^m$ and without stopping time. It yields
\begin{equation}
\begin{split}
&e^{\frac{n}{2}e^{-\alpha s}h(\| X_n^{m}(s)\|)}=e^{\frac{n}{2}
h(\|x\|)}\nonumber
\\
& \quad -\int_0^{s}\alpha \frac{n}{2} e^{-\alpha
r}h(\|X_n^{m}(r)\|) e^{\frac{n}{2}e^{-\alpha
r}h(\|X_n^{m}(r)\|)}dr\nonumber
\\
& \quad+\int_0^{s}\frac{n}{2}e^{-\alpha r}
e^{\frac{n}{2}e^{-\alpha
r}h(\|X_n^{m}(r)\|)}h'(X_n^{m}(r))\langle-A_mX_n^{m}(r)+F(X_n^{m}(r)),
\frac{X_n^{m}(r)}{\|X_n^{m}(r)\|}\rangle dr\nonumber
\\
& \quad + \int_0^{s}\frac{n}{2} e^{-\alpha r} e^{\frac{n}{2}
e^{-\alpha r}h(\|X_n^{m}(r-)\|)}h'(X_n^{m}(r-))\langle
\frac{X_n^{m}(r-)}{\|X_n^{m}(r-)\|},
G(X_n^{m}(r-))dL_n(r)\rangle\nonumber
\\
& \quad + \int_0^{s}\int_H \bigg[e^{\frac{n}{2} e^{-\alpha
r}h(\|X_n^{m}(r-)+\frac{1}{n}G(X_n^{m}(r-))z\|)}- e^{\frac{n}{2}
e^{-\alpha r}h(\|X_n^{m}(r-)\|)}\nonumber
\\
& \quad\quad -\frac{1}{2}e^{-\alpha r}e^{\frac{n}{2} e^{-\alpha
r}h(\|X_n^{m}(r-)\|)}h'(X_n^{m}(r-))\langle
\frac{X_n^{m}(r-)}{\|X_n^{m}(r-)\|},
G(X_n^{m}(r-))z\rangle\bigg]\pi_n(dr,dz).\nonumber
\end{split}\end{equation}
Arguing like in (\ref{eme221}) and (\ref{formi1}), applying
$\sup_{0\leq s\leq T}$ to both sides and taking expectation give
us

\begin{equation}\label{eme242}
\begin{split}
&\esp \sup_{0\leq s\leq T}e^{\frac{n}{2}e^{-\alpha s}h(\|
X_n^{m}(s)\|)}\leq e^{\frac{n}{2} h(\|x\|)}
\\
& \quad +\esp \sup_{0\leq s\leq T}\int_0^{s}\frac{n}{2} e^{-\alpha
r} e^{\frac{n}{2}e^{-\alpha r}h(\|X_n^{m}(r)\|)}(2C+M_2-\alpha)
h(\|X_n^{m}(r)\|) dr
\\
& \quad + \esp \sup_{0\leq s\leq T}\left|\int_0^{s}\frac{n}{2}
e^{-\alpha r} e^{\frac{n}{2} e^{-\alpha
r}h(\|X_n^{m}(r-)\|)}h'(X_n^{m}(r-))\langle
\frac{X_n^{m}(r-)}{\|X_n^{m}(r-)\|},
G(X_n^{m}(r-))dL_n(r)\rangle\right|
\\
& \quad + \esp \sup_{0\leq s\leq T}\left|\int_0^{s}\int_H
\bigg[e^{\frac{n}{2} e^{-\alpha
r}h(\|X_n^{m}(r-)+\frac{1}{n}G(X_n^{m}(r-))z\|)}- e^{\frac{n}{2}
e^{-\alpha r}h(\|X_n^{m}(r-)\|)}\right.
\\
& \left.\quad\quad -\frac{1}{2}e^{-\alpha r}e^{\frac{n}{2}
e^{-\alpha r}h(\|X_n^{m}(r-)\|)}h'(X_n^{m}(r-))\langle
\frac{X_n^{m}(r-)}{\|X_n^{m}(r-)\|},
G(X_n^{m}(r-))z\rangle\bigg]\hat\pi_n(dr,dz)\right|.
\end{split}\end{equation}

Denote
\[
N(s)=\int_0^{s}\frac{n}{2} e^{-\alpha r} e^{\frac{n}{2} e^{-\alpha
r}h(\|X_n^{m}(r-)\|)}h'(X_n^{m}(r-))\langle
\frac{X_n^{m}(r-)}{\|X_n^{m}(r-)\|}, G(X_n^{m}(r-))dL_n(r)\rangle.
\]
Then $N$ is a square integrable martingale. From the definition of
the quadratic variation process, see \cite{P},
$$
\esp[N,N]_T =\esp N^2(T).
$$
Therefore,  from the Burkholder-Davis-Gundy inequality \cite{P},
\cite{PZ},
\begin{equation}\label{eme243}
\begin{split}
& \esp \sup_{0\leq s\leq T}|N(s)|\leq C_1\esp[N,N]_T^{\frac 12}
\leq C_1(\esp[N,N]_T)^{\frac 12}=C_1(\esp N^2(T))^{\frac 12}
\\
& \leq C_2\left[\esp \int_0^{T}n^2 e^{ne^{-\alpha
r}h(\|X_n^{m}(r)\|)}\frac{M^2}{n^2}ndr\right]^{\frac 12}
\leq M_3 n^{\frac 12}
e^{\frac{n}{2} h(\|x\|)}
\end{split}\end{equation}
for some constant $M_3>0$, where we used (\ref{expest1}) to get
the last inequality. As regards the last term of (\ref{eme242}),
by Theorem 8.23 of \cite{PZ},
\begin{equation}\label{eme244}
\begin{split}
& \esp \sup_{0\leq s\leq T}\left|\int_0^{s}\int_H
\bigg[e^{\frac{n}{2} e^{-\alpha
r}h(\|X_n^{m}(r-)+\frac{1}{n}G(X_n^{m}(r-))z\|)}- e^{\frac{n}{2}
e^{-\alpha r}h(\|X_n^{m}(r-)\|)}\right.
\\
& \left.\quad\quad -\frac{1}{2}e^{-\alpha r}e^{\frac{n}{2}
e^{-\alpha r}h(\|X_n^{m}(r-)\|)}h'(X_n^{m}(r-))\langle
\frac{X_n^{m}(r-)}{\|X_n^{m}(r-)\|},
G(X_n^{m}(r-))z\rangle\bigg]\hat\pi_n(dr,dz)\right|
\\
& \quad \leq D_1 n\esp \int_0^{T}\int_H \left|e^{\frac{n}{2}
e^{-\alpha r}h(\|X_n^{m}(r)+\frac{1}{n}G(X_n^{m}(r))z\|)}-
e^{\frac{n}{2} e^{-\alpha r}h(\|X_n^{m}(r)\|)}\right.\nonumber
\\
& \left.\quad\quad\quad -\frac{1}{2}e^{-\alpha r}e^{\frac{n}{2}
e^{-\alpha r}h(\|X_n^{m}(r)\|)}h'(X_n^{m}(r))\langle
\frac{X_n^{m}(r)}{\|X_n^{m}(r)\|},
G(X_n^{m}(r))z\rangle\right|\nu(dz)dr
\\
& \leq M_4 n e^{\frac{n}{2} h(\|x\|)}
\end{split}\end{equation}
if we once again argue like in (\ref{formi1}) and then use
(\ref{expest1}).

 Therefore, plugging (\ref{eme243}) and (\ref{eme244}) in
 (\ref{eme242})we finally obtain
\begin{equation}
\esp \sup_{0\leq s\leq T}e^{\frac{n}{2}e^{-\alpha s}h(\|
X_n^{m}(s)\|)}\leq M_5 n e^{\frac{n}{2} h(\|x\|)}\leq e^{M_6 n
h(\|x\|)}
\end{equation}
for some $M_6>0$. We can now pass to the limit as $m\to+\infty$
using (\ref{yapr}) and use that $(1+r)/2\leq h(r)$ to complete the
proof.\hfill\qed

\begin{Proposition}\label{contest}
Let $0\le t\le T$ and let (\ref{as1})-(\ref{as4}) be satisfied.
Let $X_n(s)$, and $Y_n(s)$ are solutions of (\ref{a1a}) with
initial conditions $x$ and $y$ respectively. Then
\begin{equation}\label{eme11}
\esp\|X_n(s)-Y_n(s)\|^2_{-1}\le C_1(T)\|x-y\|_{-1}^2,
\end{equation}
\begin{equation}\label{eme12}
\esp\|X_n(s)-x\|^2_{-1}\le C_2(\|x\|,T)(s-t),
\end{equation}
and
\begin{equation}\label{eme13}
\esp\|X_n(s)-x\|^2\le \omega_x(s-t)
\end{equation}
for some modulus $\omega_x$.
\end{Proposition}

{\bf Proof}. The proofs are rather typical for these kinds of
estimates. We first show (\ref{eme11}). By Ito's formula we have

\begin{equation}\label{bito}
\begin{split}
& \esp\|X_n^m(s)-Y_n^m(s)\|_{-1}^2=\|x-y\|_{-1}^2
\\
& +2\esp\int_t^s[\langle X_n^m(\tau)-Y_n^m(\tau),A_m^\ast
B(X_n^m(\tau)-Y_n^m(\tau))\rangle  \\
&+\langle
F(X_n^m(\tau))-F(Y_n^m(\tau)),B(X_n^m(\tau)-Y_n^m(\tau))\rangle]d\tau
\\
& +{\frac
1n}\esp\int_t^s\int_H\|[G(X_n^m(\tau))-G(Y_n^m(\tau))]z\|_{-1}^2\nu(dz)d\tau.
\end{split}\end{equation}
Using (\ref{yapr}) and moment estimates (\ref{momest})for $X_n^m$
and $Y_n^m$ we can pass to the limit above to obtain that
(\ref{bito}) is still true if $X_n^m$ and $Y_n^m$ are replaced by
$X_n$ and $Y_n$ respectively and $A_m$ is replaced by $A$. We then
use (\ref{bcond}), (\ref{as1}) and (\ref{as2}) to get
\begin{equation}\begin{split}
& \esp\|X_n(s)-Y_n(s)\|_{-1}^2\leq \|x-y\|_{-1}^2
\\
& +(2c_0+M\|B^{\frac
12}\|)\esp\int_t^s\|X_n(\tau)-Y_n(\tau)\|_{-1}^2 d\tau
\\
& +\frac {M\|B^{\frac
12}\|}{n}\esp\int_t^s\int_H\|X_n(\tau)-Y_n(\tau)\|_{-1}^2\|z\|^2\nu(dz)d\tau
\\
& \leq \|x-y\|_{-1}^2 +C\int_t^s\esp\|X_n(\tau)-Y_n(\tau)\|_{-1}^2
d\tau \nonumber
\end{split}\end{equation}
and the claim follows from Gronwall's inequality.

To show (\ref{eme12}) we again employ Ito's formula and
(\ref{as1}), (\ref{as3}) to find that
\begin{equation}\label{bito1}
\begin{split}
& \esp\|X_n^m(s)-x\|_{-1}^2=2\esp\int_t^s[-\langle
X_n^m(\tau),A^\ast
B(X_n^m(\tau)-x)\rangle  \\
&+\langle F(X_n^m(\tau)),B(X_n^m(\tau)-x)\rangle]d\tau +{\frac
1n}\esp\int_t^s\int_H\|G(X_n^m(\tau))z\|_{-1}^2\nu(dz)d\tau
\\
& \leq C(\|x\|)\esp\int_t^s(1+\|X_n^m(\tau)\|^2)d\tau\leq
C_2(\|x\|,T)(s-t).
\end{split}\end{equation}

As regards (\ref{eme13}) it follows from the definition of mild
solution that
\[
X_n(s)=S(s-t)x+\int_t^sS(s-\tau)F(X_n(\tau))d\tau
+\int_t^sS(s-\tau)G(X_n(\tau))dL_n(\tau).
\]
Therefore
\begin{equation}\label{contmild}
\begin{split}
& \esp\|X_n(s)-x\|^2\leq
4\left[\|S(s-t)x-x\|^2+\esp\left|\int_t^sM(1+\|X_n(\tau)\|)d\tau\right|^2\right. \\
&\left.+\esp\left|\int_t^sS(s-\tau)G(X_n(\tau))dL_n(\tau)\right|^2\right]
\\
& \leq
C\left(\|S(s-t)x-x\|^2+(s-t)^2+\esp\int_t^s\frac{1}{n}d\tau\right),
\end{split}\end{equation}
where we have used the isometric formula to obtain the last
inequality. \hfill\qed

Finally we state for future use the following lemma which can be
shown rather easily using again Ito's formula applied first to the
process $X_n^m$ and then letting $m\to+\infty$. Its proof will
thus be omitted.

\begin{Lemma}\label{itof}
Let the assumptions of Proposition \ref{msol} be satisfied. Let
$t\leq s\leq T$. Let $\psi=\varphi+h(\|\cdot\|)$ be a bounded test
function. Then
\begin{equation}\begin{split}
& \esp e^{\psi(s,X_n(s))}\leq e^{\psi(t,x)}+\esp\int_t^s
e^{\psi(\tau,X_n(\tau))}[\psi_t(\tau,X_n(\tau))
\\
&+\langle F(X_n(\tau)),D\psi(\tau,X_n(\tau))\rangle +\langle
X_n(\tau),A^\ast D\varphi(\tau,X_n(\tau))\rangle] d\tau
\\
&
+n\esp\int_t^s\int_H\left[e^{\psi(\tau,X_n(\tau)+\frac{1}{n}G(X_n(\tau))z)}-e^{\psi(\tau,X(\tau))}
\right. \\
& \left. -e^{\psi(\tau,X(\tau))}\langle
D\psi(\tau,X_n\tau)),\frac{1}{n}G(X_n(\tau))z\rangle\right]\nu(dz)d\tau.
\nonumber
\end{split}\end{equation}
\end{Lemma}

\section{Associated nonlinear integro-PDE}

\label{intpde}

\setcounter{equation}{0}

For $g\in C_b(H)$ we define the function
\begin{equation} \label{bvvvf}
v_n(t,x)=\frac{1}{n}\log\esp\left(e^{ng(X_n(T))}\right),
\end{equation}
where $X_n$ solves (\ref{a1a}). As we have stated earlier one of our main aims is to establish convergence of the sequence $(v_n)$ and to
identify its limit as a solution of a Hamilton-Jacobi-Bellman equation. In the present section we investigate the approximating and the limiting
equations.

\subsection{Approximating equations}

We first show that for each $n$ the function  $v_n$ is a viscosity solution of an integro-PDE.

\begin{Theorem}\label{integropde} Let (\ref{as1})-(\ref{as4}) be satisfied and let
$g\in {\rm Lip}_b(H_{-1})$. Then there exist a constant $C_1$ and,
for every $R>0$, a constant $C_2=C_2(R)$ (both possibly depending
on $n$) such that
\begin{multline} \label{an1}
|v_n(t,x)-v_n(s,y)|\le
C_1\|x-y\|_{-1}+C_2(\max\{\|x\|,\|y\|\})|t-s|^{\frac{1}{2}}\\\mbox{for}\,\,x,y\in
H,t,s\in[0,T]
\end{multline}
and $v_n$ is a viscosity solution of  an integro-PDE
\begin{equation}
 \left\{
\begin{array}{ll}
(v_n)_t +\langle-Ax +F(x),Dv_n\rangle
& \\
& \\
\quad\quad+\int_H \left[
e^{n(v_n(t,x+\frac{1}{n}G(x)z)-v_n(t,x))}-1-\langle
Dv_n,G(x)z\rangle\right]\nu(dz) =0,
&  \\
&  \\
v_n(T,x)=g(x)\quad\quad\hbox{in}\,\, (0,T)\times H.
\end{array}
\right. \label{lali2}
\end{equation}
\end{Theorem}

{\bf Proof}. Estimate (\ref{an1}) is a direct consequence of
(\ref{eme11}), (\ref{eme12}), and the Markov property of the
process $X_n$. The proof that $v_n$ is a viscosity solution of
(\ref{lali2}) is similar to the proof of Theorem 7.1 in \cite{Sw}.
We will only show that $v_n$ is a viscosity subsolution since the
supersolution part is similar.

Suppose that $v_n-h(\|\cdot\|)-\varphi$ has a global maximum at
$(t,x)$. Since $v_n$ is bounded by Remark 4.3 of \cite{Sw} without
loss of generality we can also assume that $h, h', h''$ and
$\varphi$ are bounded. Denote $\psi(s,y)=h(\|y\|)+\varphi(s,y)$.
Then for small $\epsilon>0$
\[
v_n(t+\epsilon,X_n(t+\epsilon))- \psi(t+\epsilon,X_n(t+\epsilon))
\leq v_n(t,x)- \psi(t,x).
\]
Therefore, setting $u_n=e^{nv_n}$ we have
\[
\frac{u_n(t+\epsilon,X_n(t+\epsilon))}{u_n(t,x)} \leq e^{n
\psi(t+\epsilon,X_n (t+\epsilon))} e^{-n \psi(t,x)}.
\]
which, upon taking the expectation of both sides of the above
inequality and using the Markov property of $X_n(s)$, produces
\[
e^{n\psi(t,x)} \leq \esp e^{n \psi(t+\epsilon,X_n (t+\epsilon))}.
\]
Therefore, applying Lemma \ref{itof}, we obtain
\begin{eqnarray}
&& 0\leq \esp\frac{1}{\epsilon}\left\{ e^{n\psi(t+\epsilon,X_n
(t+\epsilon))} -e^{n\psi(t,x)}\right\} \nonumber
\\
&&
\quad\leq\esp\frac{1}{\epsilon}\int_t^{t+\epsilon}ne^{n\psi(\tau,X_n(\tau))}
\big[\psi_t(\tau,X_n(\tau))\nonumber
\\
&& \quad\quad + \langle F(X_n(\tau)),D\psi(\tau,X_n(\tau)) \rangle
d\tau - \langle X_n(\tau), A^\ast D\varphi(\tau,X_n(\tau)) \rangle
\big]d\tau \nonumber
\\
&& \quad\quad+\espe\frac{n}{\epsilon}
\int_t^{t+\epsilon}\int_H\bigg[e^{n\psi(\tau,X_n(\tau)+\frac{1}{n}G(X_n(\tau))z)}
-e^{n(\psi(\tau,X_n(\tau))} \nonumber
\\
&& \quad\quad\quad\quad -e^{n\psi(\tau,X_n(\tau))} \langle
D\psi(\tau,X_n(\tau)),G(X_n(\tau))z\rangle\bigg]\nu(dz) d\tau.
\label{ppp1gc}
\end{eqnarray}
Using (\ref{eme13}), (\ref{as1}), boundedness of $\psi$, uniform continuity of $\psi,\psi_t, D\psi, A^\ast\varphi$, and moment estimates (in particular
(\ref{eme1})) it is easy to see that
\begin{eqnarray}
&&
\esp\frac{1}{\epsilon}\int_t^{t+\epsilon}ne^{n\psi(\tau,X_n(\tau))}
\big[\psi_t(\tau,X_n(\tau))\nonumber
\\
&& \quad\quad + \langle F(X_n(\tau)),D\psi(\tau,X_n(\tau)) \rangle
d\tau - \langle X_n(\tau), A^\ast D\varphi(\tau,X_n(\tau)) \rangle
\big]d\tau \nonumber
\\
&& \quad\quad =
\frac{1}{\epsilon}\bigg[\int_t^{t+\epsilon}ne^{n\psi(t,x)}
\big[\psi_t(t,x)\nonumber
\\
&& \quad\quad \quad\quad+ \langle F(x),D\psi(t,x) \rangle d\tau -
\langle x, A^\ast D\varphi(t,x) \rangle \big]d\tau
+o(\epsilon)\bigg]. \label{rr1}
\end{eqnarray}
As regards the other term, by Lemma \ref{intform}, (\ref{as2}),
(\ref{as3}), (\ref{as4}), (\ref{eme1}), (\ref{eme13}), boundedness
of $\psi$ and uniform continuity of $\psi, D\psi, D^2\psi$, we
have
\begin{eqnarray}
&& \espe\frac{n}{\epsilon}
\int_t^{t+\epsilon}\int_H\bigg[e^{n\psi(\tau,X_n(\tau)+\frac{1}{n}G(X_n(\tau))z)}
-e^{n\psi(\tau,X_n(\tau))} \nonumber
\\
&& \quad\quad\quad\quad -e^{n\psi(\tau,X_n(\tau))} \langle
D\psi(\tau,X_n(\tau)),G(X_n(\tau))z\rangle\bigg]\nu(dz)
d\tau\nonumber
\\
&&
\quad = \espe\frac{n}{\epsilon} \int_t^{t+\epsilon}\int_H
\int_0^1\int_0^1 \langle
D^2e^{n\psi(\tau,X_n(\tau)+s\sigma\frac{1}{n}G(X_n(\tau))z)}\frac{1}{n}G(X_n(\tau))z,
\nonumber
\\
&& \quad\quad\quad\quad\quad\quad\quad\quad
\frac{1}{n}G(X_n(\tau))z\rangle\sigma
ds\,d\sigma\nu(dz)d\tau\nonumber
\\
&& \quad\quad\leq \espe\frac{n}{\epsilon}
\int_t^{t+\epsilon}\int_H \bigg[\int_0^1\int_0^1 \langle
D^2e^{n\psi(t,x+s\sigma\frac{1}{n}G(x)z)}\frac{1}{n}G(x)z,
\frac{1}{n}G(x)z\rangle\sigma ds\,d\sigma\nonumber
\\
&&
\quad\quad\quad\quad+C_1(1+\|X_n(\tau)\|^2+\|z\|^2)\|z\|^2\omega(\|X_n(\tau)-x\|(1+\|z\|))\bigg]
\nu(dz)d\tau \nonumber
\\
&& \quad\quad= \frac{n}{\epsilon}
\int_t^{t+\epsilon}\bigg[\int_H\big[e^{n\psi(t,x+\frac{1}{n}G(x)z)}
-e^{n\psi(t,x)} \nonumber
\\
&& \quad\quad\quad\quad\quad\quad\quad\quad -e^{n\psi(t,x)}
\langle D\psi(t,x),G(x)z\rangle\big]\nu(dz)+\omega_1(\epsilon)\bigg]
d\tau. \label{rr2}
\end{eqnarray}
(Above $\omega,\omega_1$ are some modului and $C_1,C_2$ are
constants, all depending on $\psi$.) Therefore plugging (\ref{rr1}) and
(\ref{rr2}) into (\ref{ppp1gc}) and sending $\epsilon\to 0$ we
obtain
\begin{eqnarray}
&& 0\leq ne^{n\psi(t,x)} \bigg(\psi_t(t,x) - \langle x,A^\ast
D\varphi(t,x)\rangle +\langle F(x),D\psi(t,x)\rangle \nonumber
\\
&& \quad\quad+ \int_H
\left[e^{n(\psi(t,x+\frac{1}{n}G(x)z)-\psi(t,x))} -1 - \langle
D\psi(t,x),G(x)z\rangle\right]\nu(dz)\bigg)\nonumber
\end{eqnarray}
which completes the proof after we divide both sides by
$ne^{n\psi(t,x)}$.\hfill \qed

\subsection{Limiting Hamilton-Jacobi-Bellman equation}

\label{controlp}

The limiting equation (obtained by letting $n\to+\infty$ in (\ref{lali2})) can be formally identified as
\begin{equation}
 \left\{
\begin{array}{ll}
v_t+\langle -Ax +F(x),Dv\rangle +H_0(G^\ast(x)Dv)=0 & \\
& \\ v(T,x)=g(x)\quad\quad\hbox{in}\,\, (0,T)\times H,
\end{array}
\right.
 \label{a3}
\end{equation}
where \[ H_0(p)=\int_H\left[e^{\langle p,z\rangle}-1-\langle
p,z\rangle\right]\nu(dz). \] It is the Bellman equation
corresponding to a deterministic control problem. For $0\leq t\leq
T$, $x\in H$, and $u(\cdot)\in M_t=\{u:[t,T]\to H: u\,\,\mbox{is
strongly measurable}\}$ we consider the state equation
\begin{equation}\label{steqdet}
X'(s)=-AX(s)+F(X(s))+G(X(s))u(s),\quad X(t)=x,
\end{equation}
and we want to maximize the cost functional
\[
J(t,x;u(\cdot))=\int_t^T-L_0(u(s))ds+g(X(T))
\]
over all controls $u(\cdot)\in M_t$, where $L_0$ is the Legendre
transform of $H_0$, i.e.
\begin{equation}\label{costl}
L_0(z)=\sup_{y\in H}\{\langle z,y\rangle-H_0(y)\}.
\end{equation}
The value function for the problem is
\begin{equation}\label{valf}
v(t,x)=\sup_{u(\cdot)\in M_t}J(t,x;u(\cdot)).
\end{equation}

The Hamiltonian $H_0$ and Lagrangian $L_0$ are both convex. By
(\ref{as4}) and the definition of $H_0$ we see that $0\leq
H_0(y)< +\infty$ for every $y\in H$, $H_0(0)=0$, and $H_0$ is
locally Lipschitz continuous on $H$. Therefore $L_0(0)=0$, $L_0(z)\geq 0$ for
every $z\in H$, and
moreover
\begin{equation}\label{growthl}
L_0(z)\geq \|z\|-H_0(\frac{z}{\|z\|})\to
+\infty\quad\mbox{as}\,\,\|z\|\to +\infty
\end{equation}
(but $L_0$ can possibly take infinite values). Since $g$ is bounded
it is then obvious that
\[
v(t,x)=\sup_{u(\cdot)\in \tilde M_t}J(t,x;u(\cdot)),
\]
where
\begin{equation}\label{mtilde}
\tilde M_t=\{u(\cdot)\in M_t: \int_t^T L_0(u(s))ds\leq
K=2\|g\|_\infty\}.
\end{equation}

We will need the following simple lemma.
\begin{Lemma}\label{llegtr}
For every $\epsilon>0$ there exists a constant $N_\epsilon
=N_\epsilon(\nu)$ such that for every $z\in H$
\[
\|z\|\le \epsilon L_0(z)+N_\epsilon.
\]
\end{Lemma}
{\bf Proof}. It follows from (\ref{costl}), (\ref{as4}), and
$L_0(0)=0$ that
\[
\|z\|=\langle \epsilon z,\frac{z}{\epsilon\|z\|}\rangle\leq
L_0(\epsilon z)+H_0(\frac{z}{\epsilon\|z\|})\leq \epsilon
L_0(z)+N_\epsilon.
\]
\hfill\qed

\begin{Lemma}\label{contestdet} Let (\ref{as1})-(\ref{as3}) be satisfied.
Let $0\le t\le T$ and $u(\cdot)\in \tilde M_t$. Then:

\smallskip\noindent
$(i)$ There exists a unique mild solution $X\in C([t,T];H)$ of
(\ref{steqdet}). Moreover there exists a constant $C_1=C_1(T,K,M)$
such that
\begin{equation}\label{eme12100}
\sup_{t\le s\le T}\|X(s)\|\le C_1(1+\|x\|).
\end{equation}

\smallskip\noindent
$(ii)$ There exists a constant
$C_2=C_2(T,K,M,c_0,\|B^{\frac{1}{2}}\|)$, such that if $X$, and
$Y$ are solutions of (\ref{steqdet}) with initial conditions $x$
and $y$ respectively then
\begin{equation}\label{eme111}
\|X(s)-Y(s)\|_{-1}\le C_2\|x-y\|_{-1}\quad\mbox{for} \,\,t\leq
s\leq T,
\end{equation}

\smallskip\noindent
$(iii)$ For every $R>0$ there exists a modulus $\omega_R$,
depending on $R,K,T,\|A^*B\|$, such that if $\|x\|\leq R$ then
\begin{equation}\label{eme121}
\|X(s)-x\|_{-1}\le \omega_R(s-t)\quad\mbox{for}\,\, t\leq s\leq T,
\end{equation}
and for every $x\in H$ there exists a modulus $\omega_x$,
independent of $u(\cdot)$, such that
\begin{equation}\label{eme131}
\|X(s)-x\|\le \omega_x(s-t)\quad\mbox{for} \,\,t\leq s\leq T.
\end{equation}
\end{Lemma}

{\bf Proof}. We first notice that by Lemma \ref{llegtr}
(applied with $\epsilon=1$)
\begin{equation}\label{l1bound}
\int_t^T\|u(\tau)\|d\tau\le K+N_1
\end{equation}
for every $u(\cdot)\in\tilde M_t$. Therefore the existence and
uniqueness of a mild solution of (\ref{steqdet}) and estimate
(\ref{eme12100}) are well known. We refer for instance to
\cite{LY}, Chapter 2, Proposition 5.3.

To show (\ref{eme111}) we notice that
\[
\begin{split}
& \|X(s)-Y(s)\|_{-1}^2=\|x-y\|_{-1}^2-2\int_t^s\langle A^*
B(X(\tau)-Y(\tau)),X(\tau)-Y(\tau)\rangle d\tau
 \\
 &
 +2\int_t^s\langle
B(X(\tau)-Y(\tau)),F(X(\tau))-F(Y(\tau))+(G(X(\tau))-
G(Y(\tau)))u(\tau)\rangle d\tau
\end{split}
\]
and therefore using (\ref{bcond}), (\ref{as1}) and (\ref{as2}) we
have
\[
\|X(s)-Y(s)\|_{-1}^2\le \|x-y\|_{-1}^2+C\int_t^s
\|X(\tau)-Y(\tau)\|_{-1}^2(1+\|u(\tau)\|)d\tau.
\]
Therefore (\ref{eme111}) follows from (\ref{l1bound}) and
Gronwall's inequality.

To prove (\ref{eme121}) we write
\[
\begin{split}
& \|X(s)-x\|_{-1}^2=-2\int_t^s\langle A^*
B(X(\tau)-x),X(\tau)\rangle d\tau
 \\
 &
 +2\int_t^s\langle
B(X(\tau)-x),F(X(\tau))+G(X(\tau))u(\tau)\rangle d\tau
\end{split}
\]
and thus using (\ref{as1})-(\ref{as3}), (\ref{eme12100}) and Lemma
\ref{llegtr} we obtain
\[
\begin{split}
& \|X(s)-x\|_{-1}^2\leq \int_t^s C_R(1+\|u(\tau)\|)d\tau
\\
& \leq \epsilon C_R \int_t^s L(u(\tau))d\tau +C_RN_\epsilon
(s-t)\leq \epsilon C_RK +C_RN_\epsilon (s-t).
\end{split}
\]
Therefore we obtain (\ref{eme121}) with
\[
\omega_R(\tau)=\inf_{\epsilon>0} (\epsilon C_RK
+C_RN_\epsilon\tau)^{\frac 12}.
\]
Estimate (\ref{eme131}) is proved similarly noticing that
\[
\|X(s)-x\|\leq\|S(s-t)x-x\|+ \int_t^s C_R(1+\|u(\tau)\|)d\tau.
\]
\hfill \qed

The definition of viscosity solution of (\ref{a3}) is the same as
Definition \ref{defvisc2-t} after we disregard the nonlocal part
and of course it is enough to have test functions which are only
once continuously differentiable. For more on viscosity solutions
of first order PDE in Hilbert spaces we refer to \cite{CL4, CL5,
LY}.

\begin{Theorem}\label{vercontpr} Let (\ref{as1})-(\ref{as3}) be satisfied and let
$g\in {\rm Lip}_b(H_{-1})$. There exist a constant $D_1$ and, for
every $R>0$, a modulus $\omega_R$ such that the value function $v$
satisfies
\begin{equation} \label{an112}
|v(t,x)-v(s,y)|\leq
D_1\|x-y\|_{-1}+\omega_R(|t-s|)\quad\mbox{for}\,\,x,y\in
H,\|x\|,\|y\|\leq R, t,s\in[0,T].
\end{equation}
Moreover $v$ is a viscosity solution of the HJB equation
(\ref{a3}).
\end{Theorem}

{\bf Proof}. The proof is very similar to the
proof of Theorem 7.3 in \cite{Sw}. We include it here for
completeness.

The Lipschitz continuity in $x$ follows from (\ref{eme111}) and
the fact that $g\in {\rm Lip}_b(H_{-1})$. To show the continuity
in time let $x\in H$ and $s<t$ and let $\epsilon>0$. Let
$u_\epsilon(\cdot)\in M_t$ be such that
\[
v(t,x)\leq J(t,x;u_\epsilon(\cdot))+\epsilon.
\]
Extending $u_\epsilon(\cdot)$ by $0$ to $[s,T]$ we can assume that
$u_\epsilon(\cdot)\in M_s$. Therefore
\[
\begin{split}
& v(s,x)-v(t,x)\geq
J(s,x;u_\epsilon(\cdot))-J(t,x;u_\epsilon(\cdot))-\epsilon
\\
& \geq g(X(T;s,x))-g(X(T;t,x))+\epsilon\geq -C_2D_2\omega_R(|s-t|)
- \epsilon,
\end{split}
\]
where we have used (\ref{eme111}), (\ref{eme121}), and $D_2$ is
the Lipschitz constant of $g$. For the opposite inequality if
$u_\epsilon(\cdot)\in M_s$ is such that
\[
v(s,x)\leq J(s,x;u_\epsilon(\cdot))+\epsilon
\]
then $u_\epsilon(\cdot)\in M_t$ and by
(\ref{eme111}), (\ref{eme121}) we again have
\[
\begin{split}
& v(s,x)-v(t,x)\leq
J(s,x;u_\epsilon(\cdot))+\epsilon-J(t,x;u_\epsilon(\cdot))
\\
& \leq
g(X(T;s,x))-g(X(T;t,x))-\int_s^tL_0(u_\epsilon(\tau))d\tau+\epsilon
\\
&
\leq C_2D_2\omega_R(|s-t|)+\epsilon. \end{split}
\]
Therefore since $\epsilon$ was arbitrary we have obtained
\[
|v(s,x)-v(t,x)|\leq C_2D_2\omega_R(|s-t|). \]

We will only show that $v$ is a viscosity subsolution as the proof
of the supersolution property is similar but easier. We will use
the dynamic programming principle. It asserts that if $0\le
t<t+\epsilon\le T, x\in H$ then
\[
v(t,x)=\sup_{u(\cdot)\in
M_t}\left\{\int_t^{t+\epsilon}-L_0(u(s))ds+v(t+\epsilon,X(t+\epsilon))
\right\}.
\]
Let now $v-\varphi-h(\|\cdot\|)$ have a local maximum at $(t,x)$.
By the dynamic programming principle for every $0<\epsilon <T-t$
there exists a control $u_\epsilon(\cdot)$ such that.
\[
v(t,x)\leq \int_t^{t+\epsilon}-L_0(u_\epsilon(s))ds+v(t+\epsilon,
X_\epsilon(t+\epsilon)) +\epsilon^2
\]
We recall that in particular this implies that $u_\epsilon(\cdot)$
is integrable.

 Denote $\psi(s,y)=-\varphi(s,y)-h(\|y\|)$. For simplicity we will write
 $h(y):=h(\|y\|)$.

We have
\[
\begin{split}
& \varphi(t+\epsilon,X_\epsilon(t+\epsilon))=
\varphi(t,X_\epsilon(t))+\int_t^{t+\epsilon} [-\langle
X_\epsilon(s),A^\ast D\varphi(X_\epsilon(s)) \rangle \\
&+\langle F(X_\epsilon(s))+
G(X_\epsilon(s))u_\epsilon(s),D\varphi(X_\epsilon(s)) \rangle] ds
\end{split}
\]
and
\[
h(X_\epsilon(t+\epsilon))\leq h(x)+\int_t^{t+\epsilon}\langle
F(X_\epsilon(s))+ G(X_\epsilon(s))u_\epsilon(s),Dh(X_\epsilon(s))
\rangle ds.
\]
The first equality above is proved for instance in \cite{LY},
Chapter 2, Proposition 5.5 and the inequality is also standard and
can be shown using Yosida approximations similarly to what
we have done in the stochastic case.

Using this we therefore have
\begin{eqnarray}
&& -\epsilon\leq \frac{1}{\epsilon}
(v(t+\epsilon,X_\epsilon(t+\epsilon))-v(t,x))-\int_t^{t+\epsilon}L_0(u_\epsilon(s))ds
\nonumber
\\
&& \leq
\frac{1}{\epsilon}(\varphi(t+\epsilon,X_\epsilon(t+\epsilon))-\varphi(t,x)
+h(X_\epsilon(t+\epsilon))-h(x))-\int_t^{t+\epsilon}
L(u_\epsilon(s))ds \nonumber
\\
&& \quad \leq\frac{1}{\epsilon}\bigg\{\int_t^{t+\epsilon}
\bigg[\varphi_t(s,X_\epsilon(s)) -\langle X_\epsilon(s), A^\ast
D\varphi(s,X_\epsilon(s)) \rangle \nonumber
\\
&& \quad\quad+\langle F(X_\epsilon(s))+
G(X_\epsilon(s))u_\epsilon(s),D\psi(s,X_\epsilon(s)) \rangle
-L_0(u_\epsilon(s))\bigg]ds\bigg\} \nonumber
\\
&& \quad \leq \frac{1}{\epsilon}\bigg\{\int_t^{t+\epsilon}
\bigg[\varphi_t(s,X_\epsilon(s)) -\langle X_\epsilon(s), A^\ast
D\varphi(s,X_\epsilon(s)) \rangle \nonumber
\\
&& \quad\quad +\langle F(X_\epsilon(s)),D\psi(s,X_\epsilon(s))
\rangle +H_0(G^\ast(X_\epsilon(s)) D\psi(s,X_\epsilon(s)))
\bigg]ds\bigg\}. \label{pppp1}
\end{eqnarray}
Therefore, using (\ref{eme131}), we can pass to the limit as
$\epsilon\to 0$ in (\ref{pppp1}) to obtain
\[
0\leq \psi_t(t,x) -\langle x, A^\ast D\varphi(t,x)\rangle+\langle
F(x),D\psi(t,x) \rangle +H_0(G^\ast(x) D\psi(t,x)).
\]
\hfill \qed

\section{Existence of Laplace limit}

\label{bcontcomparison}

\setcounter{equation}{0}

Define
\[
H(x,p)=H_0(G^*(x)p).
\]
By (\ref{as2}), (\ref{as3}) and local Lipschitz continuity of
$H_0$ we have that for every $R>0$ there exists a constant $K_R$
such that
\begin{equation}\label{Hcont}
|H(x,p)-H(y,q)|\leq K_R(\|x-y\|_{-1}+\|p-q\|)\quad\mbox{for
all}\,\,x,y,p,q\in H, \|p\|,\|q\|\leq R.
\end{equation}

The theorems below are our key results on the existence of the Laplace limit.

\begin{Theorem}\label{thcompcp1}
Let (\ref{as1})-(\ref{as4}) hold. Let $g\in {\rm Lip}_b(H_{-1})$.
Let $v_n$ be bounded viscosity solutions of (\ref{lali2}), and $v$
be a bounded viscosity solution of (\ref{a3}) such that
\begin{equation}
\lim_{t \to T}\{|v_n(t,x)-g(x)|+|v(t,x)-g(x)| \}= 0,
\textrm{uniformly on bounded sets}\label{bc1}
\end{equation}
for every $n$ and
\begin{equation}\label{lipcont2}
|v(t,x)-v(t,y)|\le D_1\|x-y\|_{-1}
\end{equation}
for some $D_1\ge 0$ and all $t\in (0,T], x,y\in H$. Let
$K:=\|v\|_\infty+\sup_n\|v_n\|_\infty <+\infty$. Then
\begin{equation}\label{supnormest}
\|v_n-v\|_\infty\to 0\quad\mbox{as}\,\,\,n\to+\infty.
\end{equation}
\end{Theorem}

The proof of this theorem is postponed until the end of the section.

\begin{Remark}\label{remuniq}
We point out that Theorem \ref{thcompcp1} implies that if (\ref{as1})-(\ref{as4}) hold and $g\in {\rm Lip}_b(H_{-1})$ then the value function
(\ref{valf}) of the control problem of Section \ref{controlp} is the unique bounded viscosity solution of (\ref{a3}) satisfying (\ref{bc1}) and
(\ref{lipcont2}).
\end{Remark}

Let $X_n(T)$ be  the solution of (\ref{a1}) (i.e. the solution of (\ref{a1a}) with $t=0$). Theorems \ref{integropde}, \ref{vercontpr}, and
\ref{thcompcp1} yield the following corollary.
\begin{Corollary}\label{corlimit}
Let (\ref{as1})-(\ref{as4}) hold and let $g\in {\rm Lip}_b(H_{-1})$. Then
\[
\Lambda(g):=\lim_{n\to\infty}\frac{1}{n}\log \esp e^{ng(X_n(T))}=v(0,x),
\]
where $v$ is the value function defined by (\ref{valf}).
\end{Corollary}
This result can now be easily extended to larger class of functions $g$.
\begin{Theorem}\label{thlimit}
Let (\ref{as1})-(\ref{as4}) hold and let $g$ be bounded and weakly sequentially continuous on $H$. Then $\Lambda(g)$ exists and
\begin{equation}\label{limg1}
\Lambda(g)=v(0,x),
\end{equation}
where $v$ is the value function defined by (\ref{valf}).
\end{Theorem}

{\bf Proof}. We use exponential moment estimate (\ref{eme1}) and the fact that $g$ can be approximated uniformly on balls in $H$ by functions in ${\rm
Lip}_b(H_{-1})$. Since (\ref{limg1}) is true for every $g\in {\rm Lip}_b(H_{-1})$, it will be preserved in the limit. Since the argument is rather
standard it will not be repeated here. Instead we refer to the proofs of Lemma 7.6 and Proposition 7.7 of \cite{Sw}. \hfill\qed

We now pass to the proof of Theorem \ref{thcompcp1}.

{\bf Proof of Theorem \ref{thcompcp1}}. If (\ref{supnormest}) is not satisfied then without loss of generality we can assume that there exists $\epsilon>0$ and a subsequence $n_k$
such that
\begin{equation}
\sup(v_{n_k}-v)\geq 4\epsilon. \label{contrconv}
\end{equation}
Let $a>0$ be such that $aT\leq \epsilon$ and let $m>0$ be such that
\[
m\geq K+\frac
{D_1^2}{\epsilon},\quad\mbox{and}\quad\frac{2D_1^2}{m}(c_0+M\|B^{\frac
12}\|) +{\frac {D_1}{m}}K_{2D_1\|B^{\frac 12}\|+1}\leq \frac{a}{2}.
\]
Let $\psi:[0,+\infty)\to [0,+\infty)$ be a smooth and nondecreasing function such that
$\psi(r)=r^2$ dor $0\leq r\leq 1$ and $\psi(r)=2$ for $r\geq 2$.
For each $k$ we choose $\mu_k>0$ such that
\[
\sup(v_{n_k}-v-\frac{\mu_k}{t}-\frac{\mu_k}{s})\geq 3\epsilon .
\]
For $\delta,\beta >0$ we now consider the function
the function
\begin {equation}
\begin{split}
\label{eqphi1} \Phi(t,s,x,y) =
v_{n_k}(t,x)&-v(s,y)-a(T-t)-\frac{\mu_k}{t}-\frac{\mu_k}{s}
-m\psi(\|x-y\|_{-1}^{2}) -{\frac{{(t-s)^{2}}}{{2\beta }}}
\\
& -\delta\sqrt{1+\|x\|^2}-\delta\sqrt{1+\|y\|^2}.
\end{split}\end{equation}
Since $\Phi$ is $B$-upper semicontinuous,

By a perturbed optimization technique of \cite{CL5}
(see page 424 there or \cite{LY}, Chapter 6.4), which is a version of the Ekeland-Lebourg
Lemma \cite{EL}, we obtain for every sufficiently big $i>0$
elements $p_i,q_i\in H$ and $a_i,b_i\in{\mathbb R}$ such that $\|p_i\|+
\|q_i\|+|a_i|+|b_i|\leq 1/i$ and such that
\[
 \Phi(t,s,x,y)+a_i t+b_i s+\langle Bp_i,x\rangle+\langle Bq_i,r\rangle
\]
has a
global maximum over $[0,T]\times H$
at some points $\bar{t},\bar{s},\bar {x},\bar{y}$, where
$0<\bar{t},\bar{s}$. Following standard arguments (see for instance
\cite{I}) is is easy to
see that
\begin{equation}
\limsup_{\delta\to 0} \limsup_{\beta \to 0}\limsup_{i\to +\infty}\,
\delta(\sqrt{1+\|\bar{x}\|^2}+\sqrt{1+\|\bar{y}\|^2})=0\quad
\hbox{for fixed} \,\,k, \label{eqxcond2ldp}
\end{equation}
\begin{equation}
\limsup_{\beta \to 0}\limsup_{i\to +\infty}\, {\frac{{(\bar{t}-\bar{s})^{2}}}{{2\beta
}}} =0\quad \hbox{for fixed}\,\,k,\delta. \label{eqxcond3ldp}
\end{equation}
Moreover it is clear that $\psi(\|\bar x-\bar y\|_{-1}^{2})=\|\bar x-\bar y\|_{-1}^{2}$
and, since
$\Phi(\bar t,\bar s,\bar x,\bar x)\le \Phi(\bar
t,\bar s,\bar x,\bar y)$, we obtain
\[
m\|\bar x-\bar y\|_{-1}^2\leq D_1\|\bar x-\bar
y\|_{-1}+\delta\sqrt{1+\|\bar x\|^2}+\langle q_i,\bar y-\bar x\rangle
\]
which, in light of (\ref{eqxcond2ldp}) and the fact that $\|\bar x\|,\|\bar y\|\leq
c_\delta$ for every $i$ for some constant $c_\delta$, implies
\begin{equation}
\limsup_{\delta \to 0}\limsup_{\beta \to 0}\limsup_{i\to +\infty}\, m \|\bar{x}
-\bar{y}\|_{-1}\le D_1. \label{boundeps11}
\end{equation}
Therefore, by (\ref{contrconv}), (\ref{eqxcond2ldp}), (\ref {eqxcond3ldp}),
(\ref{boundeps11})
and the definition of $m$, for
small $\delta,\beta$, and big $i$ we have $0<\bar{t},\bar{s}<T$.

We now use (\ref{eqphi1}) and the definition of viscosity solution to
obtain
\begin{equation}\begin{split}
&-a-a_i-\frac{\mu_k}{\bar t^2}+\frac{\bar t -\bar s}{\beta}-\<
\bar{x},A^*B(2m(\bar{x}-\bar{y})-p_i)\>\\
&+\<
F(\bar{x}),2mB(\bar{x}-\bar{y})+\frac{\delta\bar
x}{\sqrt{1+\|\bar x\|^2}}-Bp_i\> \\
& +\int_H\biggl[e^{n_km(\psi(\|\bar x+\frac{1}{n_k}G(\bar x)z-\bar y\|_{-1}^{2})
-\psi(\|\bar x-\bar y\|_{-1}^{2}))+\delta
n_k(\sqrt{1+\|\bar x+\frac{1}{n_k}G(\bar x)z\|^2}
-\sqrt{1+\|\bar x\|^2})-\langle Bp_i,G(\bar x)z\rangle}\\
& \quad\quad\quad\quad-1-\biggl\langle
2mB(\bar{x}-\bar{y})+\frac{\delta\bar x}{\sqrt{1+\|\bar
x\|^2}}-Bp_i,G(\bar x)z\biggr\rangle\biggr]\,\nu(dz) \geq 0 \label{subine12'}
\end{split}\end{equation}
and
\begin{equation}\begin{split}
&b_i+\frac{\mu_k}{\bar s^2} +\frac{\bar t -\bar s}{\beta}
-\langle\bar{y},A^*(2mB(\bar{x}-\bar{y}+q_i))\rangle +\<
F(\bar{y}),2mB(\bar{x}-\bar{y})-\frac{\delta\bar
y}{\sqrt{1+\|\bar y\|^2}}+Bq_i\>
\\
& \quad\quad\quad\quad +H\Bigl(\bar
y,2mB(\bar{x}-\bar{y})-\frac{\delta\bar y}{\sqrt{1+\|\bar y\|^2}}+Bq_i\Bigr)
\leq 0. \label{supine122}
\end{split}\end{equation}
But
\begin{equation}\begin{split}
& e^{n_k(\psi(\|\bar x+\frac{1}{n_k}G(\bar x)z-\bar y\|_{-1}^{2})
-\psi(\|\bar x-\bar y\|_{-1}^{2}))+\delta
n_k(\sqrt{1+\|\bar x+\frac{1}{n_k}G(\bar x)z\|^2}
-\sqrt{1+\|\bar x\|^2})-\langle Bp_i,G(\bar x)z\rangle}
\\
& \quad\quad\quad\quad= e^{\langle 2mB(\bar{x}-\bar{y})+\frac{\delta\bar x}{\sqrt{1+\|\bar
x\|^2}}-Bp_i,G(\bar x)z\rangle+\sigma_k(z)},
\end{split}\end{equation}
where for small $\delta,\beta$ and big $i$
\[
|\sigma_k(z)|\leq C_m\min(\|z\|,\frac{\|z\|^2}{n_k})
\]
for some constant $C_m$ independent of $k$.
Using this in (\ref{subine12'}) we therefore obtain that for
small $\delta,\beta$ and big $i$
\begin{equation}\begin{split}
&-a-a_i-\frac{\mu_k}{\bar t^2}+\frac{\bar t -\bar s}{\beta}-\langle
\bar{x},A^*B(2m(\bar{x}-\bar{y})-p_i)\rangle\\
&+\<F(\bar{x}),2mB(\bar{x}-\bar{y})+\frac{\delta\bar
x}{\sqrt{1+\|\bar x\|^2}}-Bp_i\> \\
& \quad\quad +H(\bar x, 2mB(\bar{x}-\bar{y})+\frac{\delta\bar
x}{\sqrt{1+\|\bar x\|^2}}-Bp_i) \geq -\int_{\{\|z\|\leq
1\}}\frac{\tilde C_m\|z\|^2}{n_k}\,\nu(dz)
\\
& \quad\quad\quad\quad +\int_{\{\|z\|> 1\}}e^{(2D_1\|B\|^{\frac
12}+1)M\|z\|}(e^{\sigma_k(z)}-1)\nu(dz) \geq
-\omega(k,\delta,\beta,i), \label{subine121}
\end{split}\end{equation}
where $\lim_{k\to+\infty}\limsup_{\delta\to 0}\limsup_{\beta\to 0}
\limsup_{i\to +\infty}\omega(k,\delta,\beta,i)=0$
by (\ref{as4}) and the\break Lebesgue dominated convergence theorem.

Combining(\ref{supine122}) and (\ref{subine121}) and using
(\ref{eqxcond2ldp}), (\ref{boundeps11}), (\ref{as1}), (\ref{as2}),
(\ref{as3}) we thus obtain
\begin{equation}
\begin{split}
a &\leq -2\frac{\mu_k}{T^2}+2m(c_0+M\|B^{\frac 12}\|)\|\bar x-\bar
y\|_{-1}^2+K_{2D_1\|B^{\frac
12}\|+1}\|\bar x-\bar y\|_{-1} +\omega_1(k,\delta,\beta,i)\\
& \quad\quad\quad\quad \leq \frac{2D_1^2}{m}(c_0+M\|B^{\frac
12}\|)+\frac{D_1}{m}K_{2D_1\|B^{\frac 12}\|+1} +\omega_2(k,\delta,\beta,i)
\\
& \quad\quad\quad\quad\quad\quad\quad\quad \leq
\frac{a}{2}+\omega_2(k,\beta,\delta,i),
\end{split}\end{equation}
where $\limsup_{k\to+\infty}\limsup_{\delta\to 0}\limsup_{\beta\to
0}\limsup_{i\to +\infty}\omega_j(k,\beta,\delta,i)=0$ for fixed $j=1,2$. This yields a
contradiction after we send $i\to+\infty,\beta\to 0,\delta\to 0$ and then $k\to+\infty$.

Similar argument gives us that $\lim_{n\to+\infty}\sup(v-v_n)=0$
and therefore (\ref{supnormest}) follows for some modulus
$\omega$. \hfill \qed

\section{Large deviation principle}
\label{lardevp}

\setcounter{equation}{0}

Let $V$ be a Hilbert space such that $H\subset V$ and $H\hookrightarrow V$ is compact.
We remark that on every closed ball in $H$, the topology of $V$ is
equivalent to the weak topology in $H$. We have the following large deviation result.

\begin{Theorem}\label{thldp}
Let (\ref{as1})-(\ref{as4}) hold. Let $T>0, x\in H,$ and let $X_n$
be the solutions of (\ref{a1}). Then the random variables $X_n(T)$
satisfy large deviation principle in $V$ with the rate function
\begin{equation}\label{ratef}
I(y)=\liminf_{z\to y}\inf_{u(\cdot)\in M_0}\left\{\int_0^T
L_0(u(s))ds:X\,\,\mbox{satisfies}\,\,(\ref{steqdet}),
X(0)=x,X(T)=z\right\}, \end{equation} (where the liminf above is
taken in the topology of $V$).
\end{Theorem}

{\bf Proof}. By Bryc's theorem (see for instance \cite{DPE}, Theorem 1.3.8)
to show that $X_n(T)$ satisfy
large deviation principle in $V$ it is enough to prove that
$X_n(T)$ are exponentially tight in $V$ and that for every $g\in
C_b(V)$ the Laplace limit $\Lambda(g)$ exists. Since closed balls
in $H$ are compact in $V$, exponential tightness of $X_n(T)$
follows from the exponential moment estimates (\ref{eme1}). Since
every $g\in C_b(V)$ is weakly sequentially continuous on $H$, the
Laplace limit $\Lambda(g)$ exists by Theorem \ref{thlimit}. It
remains to prove the representation formula for the rate function.
We recall that
\[
\Lambda(g)=  \sup_{u(\cdot)\in M_0}
\{\int_0^T-L_0(u(s))ds+g(X(T))\},
\]
where $X(0)=x$.

We have (see \cite{DPE}, page 27 or \cite{FK}, page 47)
\[
\begin{split}
I(y)&=\sup_{g\in C_b(V),g(y)=0}\{-\Lambda(g)\} \\
&=\sup_{g\in C_b(V),g(y)=0,g\geq 0}\inf_{u(\cdot)\in M_0}
\{\int_0^T L_0(u(s))ds+g(X(T))\}. \end{split}\] Denote the
right-hand side of (\ref{ratef}) by $I_1(y)$ and for $m>0$ define
the function
\[
g_m(z)=m\|z-y\|_V,
\]
where $\|\cdot\|_V$ is the norm in $V$. Then for $m,n\geq 1$
$$
\begin{aligned}
I(y)&\geq \inf_{u(\cdot)\in M_0} \biggl\{\int_0^T L_0(u(s))ds+g(X(T))\biggr\}
\\
& \geq \min\left\{\frac{m}{n},\inf_{u(\cdot)\in M_0}\biggl\{\int_0^T
L_0(u(s))ds: \|X(T)-y\|_V\leq\frac{1}{n}\biggr\}\right\}.
\end{aligned}
$$
Therefore, letting $m\to+\infty$ we obtain
$$
I(y)\geq\inf_{u(\cdot)\in M_0}\biggl\{\int_0^T L_0(u(s))ds:
\|X(T)-y\|_V\leq\frac{1}{n}\biggr\},
$$
which implies $I(y)\geq I_1(y)$. To show the reverse inequality,
for $g\in C_b(V)$ let $\omega_g^y$ be a modulus of continuity of
$g$ at $y$. Then for $n\geq 1$ we have
\begin{multline*}
\inf_{u(\cdot)\in M_0} \biggl\{\int_0^T L_0(u(s))ds+g(X(T))\biggr\}\\{}\leq
\inf_{u(\cdot)\in M_0}\biggl\{\int_0^T L_0(u(s))ds:
\|X(T)-y\|_V\leq\frac{1}{n}\biggr\}+\omega_g^y\Bigl({\frac 1n}\Bigr).
\end{multline*}
Taking the $\liminf_{n\to +\infty}$ in the above inequality and then
supremum over $g$ gives us $I(y)\geq I_1(y)$. \hfill\qed

\begin{Remark}\label{remratef}
Since if $\int_0^TL_0(u(s))ds\leq n$ the solution of (\ref{steqdet})
with $X(0)=x$ satisfies $\|X(T)\|\leq C_n$ for some absolute
constant $C_n$ it is clear that $I(y)=+\infty$ if $y\in V\setminus
H$.
\end{Remark}

In some cases $\liminf_{z\to y}$ can be removed from (\ref{ratef}). We present
below one such case.

\begin{Proposition}\label{remratef1}
Suppose that, in addition to the assumptions of Theorem
\ref{thldp}, there exists $p>1$ such that
\begin{equation}\label{remratefa1}
\|z\|^p\leq C(1+L_0(z))\quad\mbox{for all}\,\,\,z\in H,
\end{equation}
and that for every $x\in H$ and $K>0$ there exists a modulus
$\omega_{x,K}$ such that if $X$ satisfies (\ref{steqdet}),
$X(0)=x$, $\int_0^T\|u(s)\|^pds\leq K$, then
\begin{equation}\label{remratefa2}
\|X(s_1)-X(s_2)\|_V\leq\omega_{x,K}(|s_1-s_2|)\quad\mbox{for all}\,\,\,s_1,s_2\in[0,T].
\end{equation}
Then
\begin{equation}\label{ratef1a}
I(y)=\inf_{u(\cdot)\in M_0}\left\{\int_0^T
L_0(u(s))ds:X\,\,\mbox{satisfies}\,\,(\ref{steqdet}),
X(0)=x,X(T)=y\right\}. \end{equation}

\end{Proposition}

{\bf Proof}. To show (\ref{ratef1a}), suppose that $X_m$ satisfies
(\ref{steqdet}) with $u_m(\cdot)\in M_0$, $X_m(0)=x, X_m(T)=z_m$,
where $z_m\to y$ in $V$, and
\[
\int_0^TL_0(u(s))ds\to
\alpha\in\mathbb{R}\quad\mbox{as}\,\,\,m\to+\infty.
\]
Then by (\ref{eme12100}), (\ref{remratefa1}) and
(\ref{remratefa2}) the family $\{X_m\}$ is equibounded in $H$ and
equicontinuous in $V$ and since balls in $H$ are compact in $V$,
by the Arzela-Ascoli theorem a subsequence, still denoted by
$X_m$, converges uniformly in $C[0,T];V)$ to $Y:[0,T]\to H$ which
also satisfies (\ref{remratefa2}). Moreover we can assume that
$u_m\rightharpoonup u$ in $L^p(0,T;H)$ for some $u$. By the
definition of mild solution for $0\leq s\leq T$
\[
X_m(s)=S(s)x+\int_0^sS(s-\tau)(F(X_m(\tau))+G(X_m(\tau))u_m(\tau))
d\tau.
\]
Since the topology of $V$ on closed balls of $H$ is equivalent to
the weak topology in $H$, we have that $\sup_{0\leq \tau\leq
T}\|X_m(\tau)-Y(\tau)\|_{-1}\to 0$ as $m\to+\infty$, and thus
\begin{equation}\label{remratefa3}
\sup_{0\leq \tau\leq
T}(\|F(X_m(\tau))-F(Y(\tau))\|+\|G(X_m(\tau))-G(Y(\tau))\|)\to
0\quad \mbox{as}\,\,\,m\to+\infty.
\end{equation}
Therefore (\ref{remratefa3}), combined with $u_m\rightharpoonup u$
in $L^p(0,T;H)$, yields that for every $p\in H$
\[
\begin{split}
&
\<Y(s),p\>=\lim_{m\to+\infty}\<X_m(s),p\>
\\
&\quad =\<S(s)x+\int_0^sS(s-\tau)(F(Y(\tau))+G(Y(\tau))u(\tau))
d\tau,p\>.
\end{split}
\]
This means that $Y$ is the mild solution of (\ref{steqdet}) with
$u(\cdot)\in M_0$, $Y(0)=x, Y(T)=y$.

Since $u_m\rightharpoonup u$ in $L^p(0,T;H)$
\begin{equation}\label{remratefa4}
\sum_{i=1}^k\lambda_i^ku_{m^k_i}\to
u\quad\mbox{in}\,\,\,L^p(0,T;H)
\end{equation}
where for every $k\geq 1$, $\sum_{i=1}^k\lambda_i^k=1$ and
$\inf_{1\leq i\leq k}m^k_i\geq k$. Moreover, upon taking another
subsequence, we can assume that we have pointwise convergence in
(\ref{remratefa4}) a.e. on $[0,T]$. It now follows from Fatou's
lemma that
\[
\begin{split}
& \int_0^TL_0(u(s))ds= \int_0^T
\lim_{k\to+\infty}L_0(\sum_{i=1}^k\lambda_i^ku_{m^k_i}(s))ds
\\
& \leq
\liminf_{k\to+\infty}\int_0^TL_0(\sum_{i=1}^k\lambda_i^ku_{m^k_i}(s))ds
\leq
\liminf_{k\to+\infty}\sum_{i=1}^k\lambda_i^k\int_0^TL_0(u_{m^k_i}(s))ds=\alpha
\end{split}
\]
which completes the proof.\hfill\qed

\begin{Remark}\label{remratef1a}
Condition (\ref{remratefa2}) is satisfied for instance if
$S(\cdot)$ is a compact semigroup. We also remark that in the
above proof, (\ref{as1}) cannot be replaced by (\ref{as1strong})
even if (\ref{bcondstrong}) is satisfied.
\end{Remark}

\section{Examples of noise processes}
\label{examplesnoise} \setcounter{equation}0

We will consider two specific cases of small perturbations: compound Poisson
processes and subordinated Wiener processes. We will try to calculate the functions
\begin{eqnarray}
H_0(p)&=&\int_H \Bigl[e^{\<p,z\>}-1-\<p,z\>\Bigr]\,\nu(dz), \label{n1}\\
L_0(z)&=&\sup_{y\in H}\{\<z,y\>-H_0(y)\}. \label{n2}
\end{eqnarray}

\subsection{Compound Poisson noise}
Let $L$ be a compound Poisson process with the Gaussian jump measure $\nu = N(0,Q)$ with the trace class covariance operator $Q\ge0$, $\Tr Q<+\infty$. It
is easy to see, compare also Proposition 4.18 in \cite{PZ}, that the operator $Q$ is identical with  the covariance of $L$. It is well known, see e.g.
\cite{DPZ}, that in this specific case for each $k>0$
\begin{equation} \int_H\|z\|^2e^{k\|z\|}\,\nu(dz)<+\infty. \label{n3}
\end{equation}
To calculate the function $H_0(\cdot)$ remark that for  a random variable  $\xi$  such that $\L(\xi)=\nu$,
$$
\int \<p,z\>^2\,\nu(dz)=\esp\<p,\xi\>^2=\<Qp,p\>=\|Q^{1/2}p\|^2.
$$
Moreover, for a real valued random variable $\eta$ such that  $\L(\eta)=N(0,1)$,
$$
\esp e^{\la \eta}=e^{\frac12\la^2}, \qquad \la\in \R^1.
$$
Consequently
\begin{equation}
\int_H e^{\<p,z\>}\,\nu(dz)=\esp e^{\eta\|Q^{1/2}p\|}=e^{\frac12\<Qp, p\>}.
\end{equation}
Thus, in the present situation
\begin{equation} \label{n4} H_0(p)=e^{\frac12\<Qp,p\>}-1 = e^{\frac12 \|Q^{\frac12}p\|^{2}}-1
\end{equation}

We denote by $Q^{-1/2}$ the pseudo inverse of $Q^{1/2}$. Since $Q^{1/2}$ is self-adjoint we have an orthogonal decomposition $H=\overline{{\rm
Im}\,Q^{1/2}}\times {\rm Ker}\, Q^{1/2}$ and  we notice that $Q^{-1/2} z$ is the unique element $p_0\in \overline{{\rm Im}\,Q^{1/2}}$ such that
$Q^{-1/2}p_0=z$. For $x\in H$ will write $x=x_0+x^\perp$ to indicate the orthogonal decomposition of $x$. We have the following general result.

\begin{Proposition} \label{np0}
Assume that
$$
H_{0}(p)= h (\|Q^{\frac12}p\|),\,\,p\in H,
$$
where $Q$ is a trace class nonnegative operator and   $h$ is a convex,even function with the Legendre transform $l$. Then the Legendre transform $L_{0}$
of $H_{0}$ is of the form:
$$
L_0(z)=
\begin{cases}
l(\|Q^{-1/2}z\|), & \hbox{if }z\in \Im \Qp,\\
+\infty, & \hbox{if }z\notin \Im\Qp.
\end{cases}
$$
\end{Proposition}
{\bf Proof}. Let $z=z_0+z^\bot$. If $z^\bot\not= 0$ then
$$
L_0(z)=\sup_{p}\Bigl[\<z,p\> - h(|\Qp p\|)\Bigr] \geq \sup_{p^\bot\in {\rm Ker}\, Q^{1/2} }\<z^\perp,p^\perp\> - h(0) = +\infty.
$$
If $z=\Qp \bar p,$\,\,$\bar p \in {\overline{{\rm Im}\,Q^{1/2}}}=
H_{1}$, then
\[
 \begin{split}
L_0(z)&=\sup_p \left(\<z,p\>-h(\|\Qp p\|)\right)= \sup_p
\left[\<\bar p,\Qp p\>-h(\|\Qp p\|)\right]
\\
& = \sup_{v\in H_{1}  } \left[\<\bar p,v\>-h(\|v\|)\right]=\sup_{t\geq 0}\left[\sup_{\|v\| =t} ( \<\bar p, v  \> -h(t))\right]
\\
& =\sup_{t\geq 0}\left[\sup_{\|v\| =t} ( \<\bar p, {\frac{v}{\|v\|}}  \>t -h(t))\right] =\sup_{t\geq 0}(\|\bar p \|t - h(t))= l(\|\bar p\|)=
l(\|Q^{-1/2}z\|),
\end{split}
\]
as required.

Let now $z\in \overline{{\rm Im}\,Q^{1/2}}\setminus \Im\Qp$. When
restricted to $\overline{{\rm Im}\,Q^{1/2}}$, $Q^{1/2}$ is a
positive, self-adjoint, compact operator and $Q^{-1/2}$ exists in
the usual sense. Let $\{e_1,e_2,...\}$ be an orthonormal basis of
$\overline{{\rm Im}\,Q^{1/2}}$ composed of eigenvectors of
$Q^{1/2}$. Then $z_n=\sum_{i=1}^{n}\<z,e_i\>e_i\in {\rm Im}
\,Q^{1/2}$. Let $H_{n}$ be the linear subspace of $H$ spanned by
the vectors $\{e_1,..., e_{n}\}$ and $p = p_{n} +
p_{n}^\perp$,\,\, $z = z_n + z_{n}^\perp $,\,\, be the orthogonal
decompositions of $p$ and $z$ with respect to $H_n$ and
$H_{n}^{\perp}$. Thus
\[
\begin{split}
L_0(z) & = \sup_{p_{n} + p_{n}^\perp}\left[ \<z, p_{n} +
p_{n}^\perp \> - h(\|Q^{1/2}(p_{n} + p_{n}^\perp)\|)\right] \geq
\sup_{p_{n}}\left[ \<z, p_{n}   \> - h(\|Q^{1/2}p_{n} \|)\right]
\\
& \geq \sup_{p_n}\left[ \<z_{n} + z_{n}^{\perp}, p_{n} \> -
h(\|Q^{1/2}p_{n} \|)\right] \geq \sup_{p_n}\left[ \<z_{n} , p_{n}
\> - h(\|Q^{1/2}p_{n} \|)\right]
\\
&= \sup_{p}\left[ \<z_{n} , p \> - h(\|Q^{1/2}p
\|)\right]=l(\|Q^{-{\frac{1}{2}}}z_n\|).
\end{split}
\]
But the sequence $(\|Q^{-{\frac{1}{2}}}z_n\|)$ tends to $+\infty$ and since $l(+\infty) = +\infty$, \, $L(z) = +\infty$, as required. \hfill\qed

As a corollary we get the following proposition
\begin{Proposition} \label{np1}
Assume that $H_0$ is given by  \myref{n4}. Let $f:\R^1_+\to \R^1_+$ be the inverse function to $g(\sigma)=\sigma e^{\fp \sigma^2}$, $\sigma\ge0$. Then
$$
L_0(z)=
\begin{cases}
\left(\bigl[f(\|Q^{-1/2}z\|)\bigr]^2-1\right)e^{\fp[f(\|Q^{-1/2}z\|)]^2}
+1, & \hbox{if }z\in \Im \Qp,\\
+\infty, & \hbox{if }z\notin \Im\Qp.
\end{cases}
$$
\end{Proposition}

\begin{Remark} \rm
It is immediate that $f$ is a concave function and for every $0<a<2$ we have
$$
\sqrt{a \ln x}\le f(x)\le \sqrt{2 \ln x}, \quad \mbox{for large}\,\,\,x.
$$
\end{Remark}

\subsection{Subordinated Wiener process}

Take $L(t)=W(Z(t))$, $t\ge0$, where $W$ is a Wiener process on~$H$, say $\L(W(1))=N(0,Q_W)$ and $Z$ is a subordinator with the jump measure $\rho$ on
$[0,+\infty)$. Thus $Z$ is an increasing process starting from~$0$ and such that
\begin{gather}
\esp e^{-\la Z(t)}=e^{-t\psi(\la)}, \qquad \la\ge0,\nonumber\\
\psi(\la)=\g\la+\int_0^{+\infty}(1-e^{-\la\sigma})\,\rho(d\sigma),
\qquad \la\ge0, \label{n5}
\end{gather}
where $\g\ge0$ and $\int_0^1\sigma\rho(d\sigma)<+\infty$, $\int_1^{+\infty}
\rho(d\sigma)<+\infty$.
If $\g=1$, $\rho\equiv0$, then $Z(t)=t$, $t\ge0$ and we have $L$ identical
with the Wiener process~$W$.

\noindent We will assume that $\g=0$, find the function $H_0$ and check under
what assumptions on~$\rho$ the crucial condition  \myref{n3} is satisfied.

\noindent It is well known, see e.g. \cite{S}, \cite{PZ}, that for the L\'evy process $L$, the measure $\nu$ is of the form \beq \label{n6}
\nu=\int_0^{+\infty} N(0,tQ_W)\,\rho(dt).
\end{equation}
By direct calculations we get that the covariance operator $Q$ of $L$ is equal to,
\begin{equation}
Q = [\int_{0}^{+\infty} t \rho (dt)] Q_W  = [\esp Z(1)] Q_{W}.
\end{equation}
To simplify notation we will assume that
\begin{equation}
\esp Z(1) = 1, \,\,{\rm and}\,\,\,{\rm then}, \,\,\, Q_W =Q.
\end{equation}
Therefore,
\begin{multline*}
H_0(p)=\int_H (e^{\<\rho,z\>}-1)\,\nu(dz)
=\int_0^{+\infty}\biggl(\int_H(e^{\<\rho,z\>}-1)N(0,tQ)(dz)\biggr)\,\rho(dt)\\
{}=\int_0^{+\infty}\Bigl(e^{\fp t\<Qp,p\>}-1\Bigr)\,\rho(dt).
\end{multline*}
Thus \beq \label{n7} H_0(p)=h (\|\Qp p\|),\,\,{\rm where}\,\,\, h(u)= \int_{0}^{+\infty}(e^{{\frac{1}{2}}tu^2}-1)\rho (dt),\,\,u\geq 0,
\end{equation}
and Proposition \ref{np0} applies. An explicit formula for $L_0$ can be easily derived.

Note that
\begin{multline*}
I=\int_H \|z\|^2 e^{\k\|z\|^2}\,\nu(dz)=
\int_0^{+\infty}\rho(dt)\Bigl[\int_H \|z\|^2e^{\k \|z\|^2}N(0,tQ)(dz)\Bigr]\\
{}=\int_0^{+\infty}\rho(dt) \esp\bigl[\|W(t)\|^2 e^{\k\|W(t)\|^2}\bigr].
\end{multline*}
But $\L(W(t))=\L(\sqrt tW(1))$. Therefore
$$
I=\int_0^{+\infty} t\,\rho(dt)\bigl[\esp\|W(1)\|^2 e^{\k\sqrt t\|W(1)\|}\bigr].
$$
We will need the following lemma.
\begin{Lemma}
There exists $a>0$ such that for all $s\ge0$,
$$
\esp e^{s\|W(1)\|}\le e^{as^2}.
$$
\end{Lemma}
{\bf Proof}. By \cite{KW}, page 55, there exists $\d>0$ such that
$$
\P(\|W(1)\|>u)\le e^{-\d u^2},\ u>0.
$$
Therefore
$$
\esp(s^{s\|W(1)\|})=\int_0^{+\infty}\P(e^{s\|W(1)\|}\ge u)\,du
=1+\int_1^{+\infty}\P\Bigl(\|W(1)\|>\frac{\ln u}s\Bigr)\,du.
$$
Note that
$$
\int_1^{+\infty}\P\Bigl(\|W(1)\|>\frac{\ln u}s\Bigr)\,du \le
\int_1^{+\infty}e^{-\d(\ln u/s)^2}\,du.
$$
Substituting $v=\frac{\ln u}s$, $du=us\,dv=se^{vs}\,dv$,
\begin{multline*}
\int_1^{+\infty}e^{-\d(\ln u/s)^2}\,du=s\ioi e^{-\d v^2}e^{vs}\,dv
=s\biggl(\ioi e^{-\d(v-s/(2\d))^2}\,dv\biggr)e^{s^2/(4\d)}\\
{}\le s\biggl(\int_{-\infty}^{+\infty} e^{-\d v^2}\,dv\biggr) e^{s^2/(4\d)}.
\end{multline*}
The required result now follows. \hfill\qed
\begin{Proposition}
If
$$
\int_0^{+\infty} t\,\rho(dt) =1  \quad\hbox{and}\quad \int_1^{+\infty}e^{\la t}\,\rho(dt)<+\infty, \quad \la\ge0,
$$
then the measure $\nu$ given by \myref{n6} satisfies \myref{n3} and $H_0$ is given by \myref{n7}.
\end{Proposition}
{\bf Proof.} It is enough to remark that,
$$
\esp\|W(1)\|^2e^{\k\sqrt t\|W(1)\|} \le \bigl(\esp\|W(1)\|^2\bigr)^{1/2}
\bigl(\esp e^{2\k\sqrt t\|W(1)\|}\bigr)^{1/2} \le ce^{\frac a2 \k^2t}.
$$
\hfill\qed

\begin{Example} \label{ne2} \rm
The assumptions of the above proposition are satisfied if, for instance,
$$
\rho(dt)=\frac1{t^{1+\alpha}}\,e^{-t^2}\,dt \quad \hbox{for }\alpha<1.
$$
\end{Example}
In some cases asymptotic behavior of the function $\psi$ can be determined.
\begin{Example} \label{ne3}\rm
\begin{gather*}
\rho(dt)=\I_{[0,1]}(t)\,\frac1{t^{1+\alpha}}\,dt, \quad \alpha<1,\\
-\psi(-\la)=\int_0^1 (e^{\la\sigma}-1)\frac1{\sigma^{1+\alpha}}\,d\sigma.
\end{gather*}
After substitution, $\la\sigma=u$, for $\la>1$,
\begin{multline*}
\int_0^1(e^{\la\sigma}-1)\frac1{\sigma^{1+\alpha}}\,d\sigma
=\frac1\la \int_0^\la (e^u-1)\frac1{(\frac u\la)^{1+\alpha}}\,du
=\la^\alpha \int_0^\la (e^u-1)\frac1{u^{1+\alpha}}\,du\\
{}\le\la^\alpha\biggl[\int_0^1 \frac{e^u-1}u\cdot \frac1{u^\alpha}\,du
+\int_1^\la e^u\,du\biggr].
\end{multline*}
Thus, for large $\la$,
$$
\int_0^1 (e^{\la\sigma}-1)\frac1{\sigma^{1+\alpha}}\,d\sigma \sim
c\la^\alpha e^\la
$$
\end{Example}

\begin{Remark} \rm
In the considered examples, the Legendre transforms $L_0$ of $H_0$
were of the form $l(\|Q^{-{\frac{1}{2}}}z\|),\,z\in H$. Thus the
control system, which defines the rate function, can be written in
a more convenient way,
\begin{equation}
X'(s)=-AX(s)+F(X(s))+G(X(s))Q^{1/2}u(s),\quad X(t)=x,
\end{equation}
and to find the rate function one has to look for the infimum of the  cost functional
\[
J(x;u(\cdot))=\int_0^T l(u(s))ds+g(X(T))
\]
over all controls $u(\cdot)\in M_0$.
\end{Remark}

\section{Stochastic PDE of hyperbolic type} \label{wave}
\setcounter{equation}0

We present an example of a class of stochastic PDE which can be handled by the
developed theory. To begin consider a nonlinear stochastic wave equation which
can be formally written as

\begin{equation}\label{w1}
\left\{
\begin{array}{ll}
&\frac{\partial^2u}{\partial t^2}(t,\xi)=\D u(t,\xi)+f(u(t,\xi))+\der{}t\tilde L_n(t,\xi),\quad t>0,\ \xi\in
\O,\\
&u(t,\xi)=0, \quad\quad\quad\quad\quad\quad\quad\quad\quad\quad\quad\quad\quad\quad
\quad\,\,\,\,
t>0,\ \xi\in\partial \O,\\
&u(0,\xi)=u_0(\xi), \quad\quad\quad\quad\quad\quad\quad\quad\quad\quad
\quad\quad\quad\,\,\,\,\,\xi\in \O,\\
&\frac{\partial u}{\partial t}(0,\xi)=v_0(\xi), \quad\quad\quad\quad\quad\quad\quad
\quad\quad\quad\quad\quad\quad\,\,\,
\xi\in \O,
\end{array}
\right.
\end{equation}
with $\tilde L_n$, $L^2(\O)$ valued L\'evy process (properly normalized), $\O$ a
bounded regular domain in $\R^d$, $f:\R\to\R$ is a Lipschitz function and
$u_0\in H_0^1(\O)$, $v_0\in L^2(\O)$.

\noindent Setting
$$
X(t)=\binom{u(t)}{v(t)},\,\,t\ge0,
$$
we can rewrite \myref{w1} in an abstract way:
\begin{equation}\label{w2}
dX(t)=\left(\left(\begin{matrix} 0&\ &I\\-A&&0 \end{matrix}\right)
X(t)+F(X(t))\right)\,dt +dL_n(t),
\end{equation}
where
\begin{equation}\label{w3}
F\binom{u}{v}=\binom{0}{F_1(u)}, \qquad L_n(t)=\binom0{\tilde L_n(t)}
\end{equation}
and $A=-\Delta$ in $H=L^2(\O)$ with $D(A)=H^2(\O)\cap H^1_0(\O)$.
Moreover the same setup applies to other equations of hyperbolic
type.

Therefore let us assume that $A$ in (\ref{w2}) is a strictly
positive, self-adjoint operator in a Hilbert space $H$ with a
bounded inverse. It is then well known that the operator
$$
\A=\pmx{0&\ &-I\\A&&0}, \qquad \cD(\A)=\left(\mx{D(A)\\\times\\D(A^{1/2})}\right)
$$
is maximal monotone in the
Hilbert space $\H=\left(\mx{D(A^{1/2})\\\times\\ H}\right)$, equipped with the following
``energy'' type inner product
$$
\<\pmx{u\\v},\pmx{\bar u\\\bar v}\>_\H
=\<A^{1/2}u,A^{1/2}\bar u\>_H+\<v,\bar v\>_H, \qquad
\pmx{u\\v},\pmx{\bar u\\\bar v}\in \H.
$$
Moreover, $\A^\ast=- \A$.

It is easy to check that the operator
$$
\B=\pmx{A^{-1/2}&\ &0\\0&& A^{-1/2}}
$$
is bounded, positive, self-adjoint on $\H$, and such that $\A^\ast \B$ is bounded. Moreover
\myref{bcond} holds with constant $c_0=1$. In fact
\[
\<(\A^\ast +I)\B\pmx{u\\v},\pmx{u\\v}\>_\H
=\<\B\pmx{u\\v},\pmx{u\\v}\>_\H=\|A^{1/4}u\|^2+\|A^{-1/4}v\|^2.
\]
In particular we see that
$$
\left\|\pmx{u\\v}\right\|_{-1}=\left(\|A^{1/4}u\|^2+\|A^{-1/4}v\|^2
\right)^{1/2}.
$$
Thus $F=\binom 0{F_1}$ is Lipschitz from $\H_{-1}$ into $\H$ (condition
(\ref{as1})) if and only if
\begin{equation} \label{w5}
\left\|A^{-1/4}(F_1(u)-F_1(\bar u)\right\|_H \le c\|A^{1/2}(u-\bar u)\|,
\qquad u,\bar u\in D(A^{1/2}).
\end{equation}
It is easy to see that if
$$
F_1(u)(\xi)=f(u(\xi)), \qquad \xi\in \O,
$$
and $f$ is a Lipschitz function, then \myref{w5} is satisfied.

\section{Appendix: Proof of Proposition \ref{apt1}}

\setcounter{equation}0

Let us recall that  the spaces $\cal X$, $\cal L$   were
introduced  in Section $4$.  Define, for each $\psi \in \cal L$, processes
$$
\begin{aligned}
\K(\psi)(t)&=\int_0^t S(t-s)\psi(s)\,dL(s), &\qquad& t\in[0,T],\\
\K_\la(\psi)(t)&=\int_0^t S_\la(t-s)\psi(s)\,dL(s),
\quad \la>0, &\qquad& t\in[0,T].\\
\end{aligned}
$$
We can treat $\K$ and $\K_\la$ as linear transformations from the space $\cal L$ into $\cal X$. We prove this now and establish that there exists a
constant $C_1>0$ such that
\begin{equation}\label{ap3}
\|\K_\la\|\le C_1  \quad\hbox{for }\la>1.
\end{equation}
In the proof we  omit the subscript $\lambda$. Let $\wh H$, and the unitary semigroup $\wh S$, be the extensions, respectively of $H$ and of the
semigroup $S$, given by the delation theorem, see e.g. \cite[Theorem 9.24]{PZ}. Thus $H\hookrightarrow \wh H$ is an isometry and the semigroup $S$ is the
restriction of $P \wh S$ to $H$, where $P$ is the orthogonal projection of $\wh H$ onto $H$. Therefore we have:
$$
\int_0^t S(t-s)\psi(s)\,dL(s) = \int_0^t P\wh S(t-s)\psi(s)\,dL(s)= P \wh S(t) \int_0^t \wh S(-s)\psi(s)\,dL(s),\quad t\in[0,T].
$$
Moreover the process
$$
\wh Y(t)=\int_0^t \wh S(-s)\psi(s)\,dL(s), \quad t\ge0,
$$
is a $\wh H$ martingale and therefore has \ca\ modification. This implies that the stochastic convolution has $H$-valued, \ca\ modifications
and
$$
\biggl\| \int_0^t S(t-s)\psi(s)\,dL(s)\biggr\| \le \bigl\|\wh Y(t)\bigr\|_{\wh H}, \quad t\in[0,T].
$$
However, $\|\wh Y(t)\|_{\wh H}$, $ t \in[0,T]$, is a submartingale and by the classical Doob inequality for all $p>1$
$$
\esp\Bigl(\sup_{0\le t\le T}\bigl\|\wh Y(t)\bigr\|^p_{\wh H}\Bigr) \le \Bigl(\frac p{p-1}\Bigr)^p
\esp\bigl\|\wh Y(T)\bigr\|_{\wh H}^p.
$$
In particular \bmlg \esp\biggl(\sup_{0\le t\le T}\biggl\|\int_0^t S(t-s)\psi(s)\,dL(s)\biggr\|^2_H\biggr)
\le 4\esp\biggl\|\int_0^T \wh S(-s)\psi(s)\,dL(s)\biggr\|_{\wh H}^2\\
{}\le 4\esp\int_0^T \bigl\|\wh S(-s)\psi(s)\Qp\bigr\|^2_{L_{\rm HS}(H,\hat H)}\,ds \le 4\esp \int_0^T \bigl\|\psi(s)\Qp\bigr\|^2_{\rm HS}\,ds.
\end{multline*}

\noindent Thus the existence of the constant $C_1$ follows, and by the Banach-Steinhaus
theorem it is enough  to establish \er{ap1} for a dense set of
$\psi$.
\begin{Lemma}\label{apl1}
For each $k=1,2,\dots$ the set
$$
{\cal L}_k=\biggl\{\psi\in {\cal L}: \esp\int_0^T
\bigl\|A^k\psi(u)\Qp\bigr\|^2_{\rm HS} \,du<+\infty\biggr\}
$$
is dense in $\cal L $.
\end{Lemma}

{\bf Proof}. Let $\psi\in {\cal L}$. Since for $\mu> 0$ the operator $\mu AR_\mu$ is bounded we have
\[
\esp\int_0^T \bigl\|A^k(\mu R_\mu)^k\psi(u)\Qp\bigr\|^2_{\rm
HS}\,du =\esp\int_0^T \bigl\|(\mu
AR_\mu)^k\psi(u)\Qp\bigr\|^2_{\rm HS}\,du <+\infty,
\]
and thus $(\mu R_\mu)^k\psi\in {\cal L}_k$. Moreover it follows
from (\ref{resolv}) that
\[ \bigl\|\bigl((\mu
R_\mu)^k-I\bigr)\psi(u)\Qp\bigr\|^2_{\rm HS} \le C
\bigl\|\psi(u)\Qp\bigr\|^2_{\rm HS}.
\]
and $\lim_{\mu\to+\infty}(\mu R_\mu)^k x=x$ for every $x\in H$.
Therefore the dominated convergence theorem yields
\[
\lim_{\mu\to+\infty}\esp\int_0^T \bigl\|\big((\mu
R_\mu)^k-I\big)\psi(u)\Qp\bigr\|^2_{\rm HS} \,du=0.
\]
\hfill\qed

\begin{Lemma}\label{apl2}
Assume that $M(t)$, $t\ge0$, is a $D(A)$-valued process with
locally bounded trajectories, $H$-square integrable martingale,
and $M(0)=0$. Then \begin{equation}\label{ap4} \int_0^t
S(t-s)\,dM(s)=M(t)-\int_0^t S(t-s)AM(s)\,ds.
\end{equation}
\end{Lemma}

{\bf Proof}. Let $e\in D((A^*)^2)$ and
\[
\varphi(s,x)=\langle S(t-s)x,e\rangle=\langle x,S^*(t-s)e\rangle.
\]
Then $\varphi\in C^{2}((-\infty,t)\times H)$ and has uniformly
continuous derivatives. In fact it can be extended to a function
in $C^{2}(\mathbb{R}\times H)$ in an obvious way. Therefore,
applying Ito's formula for Hilbert space valued semimartingales
(see \cite[Theorem 27.2]{M} or \cite[Theorem D2]{PZ}) we obtain
\[
\langle M(t),e\rangle=\int_0^t \langle S(t-s)AM(s),e\rangle ds +
\int_0^t \langle S(t-s)dM(s),e\rangle ds
\]
which proves the claim since $D((A^*)^2)$ is dense in $H$.
\hfill\qed

Applying Lemma \ref{apl2} to the martingale $M(t)= \int_0^t \psi(u)\,dL(u)$, $t\in[0,T]$
we arrive at the following lemma.
\begin{Lemma}
If $\esp\int_0^T \|A\psi(u)\Qp\|^2_{\rm HS}\,du<+\infty$ then for
all $t\in[0,T]$, $\la>0$,
$$
\aligned \int_0^t S(t-s)\psi(s)\,dL(s)&=\int_0^t \psi(s)\,dL(s)
-\int_0^t S(t-s)\biggl(\int_0^s A\psi(u)\,dL(u)\biggr)ds,\\
\int_0^t S_\la(t-s)\psi(s)\,dL(s)&=\int_0^t \psi(s)\,dL(s)
-\int_0^t S_\la(t-s)\biggl(\int_0^s A_\la\psi(u)\,dL(u)\biggr)ds.
\endaligned
$$
\end{Lemma}

\medskip
We can now continue the proof of the theorem. We will show that
(\ref{ap1}) holds for every $\psi\in {\cal L}_2$. Note that \bmlg
\K\psi(t)-\K_\la\psi(t)=\int_0^t S(t-s)
\biggl[\int_0^s -A\psi(u)\,dL(u)+\int_0^sA_\la\psi(u)\,dL(u)\biggr]ds\\
{}+\int_0^t\bigl[S(t-s)-S_\la(t-s)\bigr]\biggl(\int_0^s -A_\la\psi(u)\,dL(u)
\biggr)ds=I^1_\la\psi(t)+I^2_\la\psi(t).
\end{multline*}
Now
$$
I^1_\la\psi(t)=\int_0^t S(t-s)(A_\lambda
-A)\int_0^s\psi(u)\,dL(u)\,ds,
$$
and
$$
\sup_{0\le t\le T}\bigl\|I^1_\la\psi(t)\bigr\|\le \int_0^T \biggl\|(A-A_\la)\int_0^s \psi(u)\,dL(u)\biggr\|ds.
$$
But $\|(A-A_\la)x\|=\|R_\la A^2x\|\le \frac1\la \|A^2x\|$, $x\in
D(A^2)$. Therefore, since \beq\label{ap9} \esp\int_0^T
\bigl\|A^2\psi(u)\Qp\bigr\|^2_{\rm HS}\,du<+\infty,
\end{equation}
we have, by isometric identity,
$$
\aligned \esp\sup_{0\le t\le T}\bigl|I^1_\la\psi(t)\bigr|^2
&\le\esp\biggl(\int_0^T\biggl|(A-A_\la)\int_0^s \psi(u)\,dL(u)\biggr|\,ds
\biggr)^2\\
&\le T\int_0^T \esp\int_0^s \bigl\|(A-A_\la)\psi(u)\Qp\bigr\|^2_{\rm HS}\,du
\,ds\\
&\le\frac1{\la^2}\,T \int_0^T \esp\int_0^s\bigl\|A^2\psi(u)\Qp\bigr\|^2_{\rm HS}
\,du\,ds\\
&\le\frac1{\la^2}\,T^2 \int_0^T \esp\bigl\|A^2\psi(u)\Qp\bigr\|^2_{\rm HS}\,du.
\endaligned
$$
Therefore, if \er{ap9} holds,
$$
\lim_{\la\to+\infty}\esp\bigl\|I^1_\la\psi(t)\bigr\|^2=0.
$$
Since for every $x\in D(A)$, $\la>0$,
\[
\|S_\la(t)x-S(t)x\|\le t\|A_\la x-Ax\|
\]
(see for instance \cite{Pazy}, page 10), we have

Thus
\begin{multline*} \sup_{0\le t\le T} \bigl\|I^2_\la\psi(t)\bigr\|^2
\le \sup_{0\le t\le T}\biggl(\int_0^t\biggl\|\bigl[S(t-s)-S_\la(t-s)\bigr]A_\la
\int_0^s \psi(u)\,dL(u)\biggr\|\,ds\biggr)^2\\
{}\le\sup_{0\le t\le T}\biggl[\int_0^t (t-s)\bigl\|(A-A_\la)A_\la
\int_0^s\psi(u)\,dL(u)\bigr\|\,ds\biggr]^2\\
{}\le T^2\sup_{0\le t\le T}\biggl[\int_0^t \bigl\|(A-A_\la)A_\la
\int_0^s\psi(u)\,dL(u)\bigr\|\,ds\biggr]^2.
\end{multline*}
Moreover,
$$
(A-A_\la)A_\la=(A-\la R_\la A)\la R_\la A=\la R_\la(I-\la R_\la)
A^2.
$$
Therefore
\begin{multline*} \esp\sup_{0\le t\le T}\bigl\|I^2_\la \psi(t)\bigr\|^2
\le T^2\esp \int_0^T \biggl\|(I-\la R_\la)A^2 \int_0^s \psi(u)\,dL(u)
\biggr\|^2\,ds\\
{}\le T^2\esp\int_0^T \int_0^s \bigl\|(I-\la R_\la)A^2\psi(u)\Qp\bigr\|^2
_{\rm HS}\,ds\,du\\
{}\le T^3\esp \int_0^T \bigl\|(I-\la R_\la)A^2\psi(u)\Qp\bigr\|^2 _{\rm HS}\,du.
\end{multline*}
Thus, if \er{ap9} holds, we can conclude by the dominated
convergence theorem that
$$
\lim_{\la\to+\infty}\esp\sup_{0\le t\le T}\bigl\|I^2_\la \psi(t)\bigr\|^2=0.
$$
This finishes the proof of the proposition. \hfill\qed

\end{document}